\documentclass[11pt,table,reqno]{amsart}

\usepackage{amsfonts,amssymb,graphicx,color,tikz,amsmath,amsthm}

\usepackage{todonotes}
\usepackage{hyperref}

\usepackage{subcaption}

\tikzset{vertex/.style={circle,draw,fill,inner sep=0pt,minimum size=1mm}}
\tikzset{comp-vertex/.style={vertex,fill=white}}

\tikzset{pixel/.style={draw=black, fill=lightgray, shape=rectangle, minimum size=19 }}

\tikzset{point/.style={circle,draw,fill,inner sep=0pt,minimum size=0.25mm}} 

\tikzset{edge/.style={}}
\tikzset{pixelmode/.append style={scale=1}}

\tikzset{every picture/.append style={scale=.5}}

\usetikzlibrary{arrows.meta}

\newtheorem{thm}{Theorem}

\newtheorem{lem}[thm]{Lemma}
\newtheorem{prop}[thm]{Proposition}
\newtheorem{conj}[thm]{Conjecture}
\newtheorem{quest}[thm]{Question}

\theoremstyle{definition}
\newtheorem{definition}[thm]{Definition}
\newtheorem{exa}[thm]{Example}

\numberwithin{thm}{section}

\newcommand{\Z}{\mathbb{Z}}

\newcommand{\adj}{\sim}

\DeclareMathOperator{\id}{id}

\usepackage{tabularx}

\newcommand{\vcbox}[1]{\vcenter{\hbox{#1}}}

\begin{document}

\title{Homotopy equivalence of digital pictures in $\Z^2$}

\author{Dae-Woong Lee, and P. Christopher Staecker}

\date{\today}

\address{
Department of Mathematics, and Institute of Pure and Applied Mathematics, Jeonbuk National University,
567 Baekje-daero, Deokjin-gu, Jeonju-si, Jeollabuk-do 54896, Republic of Korea
}
\email{dwlee@jbnu.ac.kr}

\address{
Mathematics Department, Fairfield University, 1703 North Benson Rd, Fairfield, CT 06824-5195
}
\email{cstaecker@fairfield.edu}


\begin{abstract}
We investigate the properties of digital homotopy in the context of digital pictures $(X,\kappa,\bar \kappa)$, where $X\subsetneq \Z^n$ is a finite set, $\kappa$ is an adjacency relation on $X$, and $\bar \kappa$ is an adjacency relation on the complement of $X$. In particular we focus on homotopy equivalence between digital pictures in $\Z^2$.

We define a numerical homotopy-type invariant for digital pictures in $\Z^2$ called the outer perimeter, which is a basic tool for distinguishing homotopy types of digital pictures.
When a digital picture has no holes, we show that it is homotopy equivalent to its rc-convex hull, obtained by ``filling in the gaps'' of any row or column. We show that a digital picture $(X,c_i,c_j)$ is homotopy equivalent to only finitely many other digital pictures $(Y,c_i,c_j)$. At the end of the paper, we raise a conjecture on the largest digital picture of the same homotopy-type of a given digital picture.
\end{abstract}

\maketitle

\section{Introduction}

Digital topology began in the late 1960s in the work of A. Rosenfeld and others. It mainly deals with the topological properties of digital or discrete spaces, which consist of points with integer coordinates with suitable adjacency relations in $\mathbb{Z}^n$. Digital topology adapts concepts like connectedness, surfaces, curves, manifolds, holes, boundaries, homotopy, (co)homology, etc. to discrete settings, and it often uses tools from graph theory, algebraic topology, and combinatorial topology. The fundamental concepts and the major theoretical results of digital topology through the 1980s are surveyed by T. Y. Kong and A. Rosenfeld \cite{kr89}; see also \cite{krr92}.



At a basic level, Rosenfeld's digital topology views a digital image as a finite set of pixels together with some relation describing which pixels are adjacent to each other. These adjacencies are typically based on the standard integer lattice, making some choices about which diagonal adjacencies to allow.

Viewing the pixels as vertices and the adjacencies as edges, one can view a digital image as a finite graph. One of the fundamental ideas in the digital topology of the 1970s and 80s was the counterintuitive decision to consider both the digital image $X$ and its complement $\Z^2 - X$, and use different adjacency relations on each of these.

Digital homotopy between digitally continuous functions was introduced by L. Boxer \cite{boxe94}.
Work based on this notion has continued to the present, and developed in parallel with $A$-theory \cite{bl05} and $\times$-homotopy theory \cite{doch09b}. These two theories have defined invariants and concepts in the context of abstract graph theory that in many cases match ideas independently developed for the topology of digital images.


This work in graph theory was, until recently, not known in the digital topology literature, and in many instances was addressing the same questions at a more sophisticated level than had been done in digital topology. So there is now a need in the digital topology field to emphasize ideas which are truly distinctive of our point of view. This paper aims to do this by returning to the idea of complements and complementary adjacency relations. This idea is completely absent in the graph theory literature, because an abstract graph $G$ has no well-defined complement, while a digital image $X\subset \Z^2$ always does.

In this paper we ask and answer some natural questions about homotopy equivalence in the sense of Boxer (and $A$-theory and $\times$-homotopy theory) in the context of digital images with their complements using complementary adjacencies. We focus in this paper on digital pictures in $\Z^2$. Along with obvious questions about generalizing to $\Z^n$ for $n>2$, we pose many questions for further research.

This paper is organized as follows. In Section \ref{basicdef}, we describe some material on digital images in $\Z^2$ with $c_1$ and $c_2$ adjacencies, and the basic properties of digital $1$- and $2$-homotopy.
In Section \ref{complementsandpicture}, we define a digital picture as a triple $(X,\kappa,\bar \kappa)$, where $X\subsetneq \Z^2$ is a finite set, $\kappa$ is an adjacency relation on $X$, and $\bar \kappa$ is an adjacency relation on the complement $\bar X$ of $X$. We also define isomorphism and homotopy equivalence between digital pictures.
In Section \ref{digitalrectangles}, we define the outer perimeter, which we prove is a numerical homotopy-type invariant for digital pictures. We discuss the isomorphism type of rectangles, and show that homotopy equivalent rectangles have the same perimeter.
In Section \ref{rowcomponent}, we consider the row component graph, which is obtained by viewing a digital image $(X,c_i)$ as a graph and taking the quotient which identifies any two points in the same row component.
In Section \ref{rcconvex}, we define row-convexity, column-convexity, and the row-convex hull of a digital picture. We show that a digital picture with no hole is homotopy equivalent to its row-column-convex hull.
In Section \ref{homotopytypeofdigitalpicture}, we revisit the concepts of $i$-reducibility in digital homotopy theory, and show that a given digital picture $(X,c_i,c_j)$ is homotopy equivalent to only finitely many other digital pictures $(Y,c_i,c_j)$. We raise a question about $1$-reducible digital images, and a conjecture on the row-column-convex hull of a digital picture.

\bigskip

\section{Basic definitions} \label{basicdef}

A \emph{digital image} $(X, \kappa)$ consists of a subset $X$ of points in $\Z^n$ with some \emph{adjacency relation} $\kappa$, which is reflexive and symmetric. The typical adjacencies used in $\Z^n$ are the $c_i$-adjacencies, defined as follows:

\begin{definition}(\cite{boxe06})
Two elements $x =(x_1,x_2,\ldots, x_n)$ and $y =(y_1,y_2,\ldots,y_n)$ of $\mathbb Z^n$ are \emph{$c_u$-adjacent} if
\begin{itemize}
\item there are at most $u$ distinct indices $i$ such that $\vert x_i - y_i\vert =1$; and
\item $x_j = y_j$ for all other indices $j$.
\end{itemize}
\end{definition}

In this paper we focus on digital images $X\subset \Z^2$, for which there are two natural adjacencies: $c_1$-adjacency, in which each point is adjacent to its 4 horizontal and vertical neighbors, and $c_2$-adjacency, in which each point is adjacent to its 8 neighbors (including diagonals). In the classical literature these are called 4-adjacency and 8-adjacency, but as we will see the $c_i$ notation is more convenient for our results.

In a digital image $(X,\kappa)$, we write $x\sim_\kappa y$ when $x,y\in X$ are $\kappa$-adjacent. When the adjacency relation is clear from context we simply write $x\sim y$.

A digital image $(X,\kappa)$ is equivalent to a reflexive graph with vertex set $X$ and an edge connecting $x,y\in X$ when $x\sim y$. We will sometimes draw a digital image as a graph to emphasize the specific adjacencies, or sometimes simply as a collection of square pixels in $\Z^2$.

We will use various simple concepts from graph theory, in particular we say $(X,\kappa)$ is \emph{$\kappa$-connected} when it is connected as a graph, and we will discuss components, paths, and degree of a vertex all using their usual meaning in the graph theoretic setting.

\begin{definition}(\cite{boxe99})
A function $f : (X,\kappa_X) \rightarrow (Y,\kappa_Y)$ between digital images is said to be a [digitally] \emph{$(\kappa_X,\kappa_Y)$-continuous function} when: if $x_1 \adj x_2$ in $X$, then  $f(x_1) \sim f(x_2)$ in $Y$ for any $x_1, x_2\in X$.
\end{definition}

If the adjacency relations are clear, then we simply say that $f$ is [digitally] \emph{continuous}. It can be seen that any composition of digitally continuous functions is also digitally continuous. When $X$ and $Y$ are viewed as reflexive graphs, digital continuity of $f$ is equivalent to $f$ being a homomorphism of graphs.

A $(\kappa_X,\kappa_Y)$-continuous function $f : (X,\kappa_X) \rightarrow (Y,\kappa_Y)$ is called a \emph{$(\kappa_X,\kappa_Y)$-isomorphism} if $f$ is a bijection, and its inverse $f^{-1} : (Y,\kappa_Y) \rightarrow (X,\kappa_X)$ is $(\kappa_Y,\kappa_X)$-continuous. In this case, the two digital images $(X,\kappa_X)$ and $(Y,\kappa_Y)$ are said to be \emph{$(\kappa_X,\kappa_Y)$-isomorphic}; see  \cite{boxe94} and \cite{boxe05}.


Let $a$ and $b$ be integers with $a \le b$. A {\it digital interval} \cite {boxe94} is a finite subset of $\Z$ of the following type:
$$
[a,b]_{\mathbb{Z}} = \{ x \in \mathbb Z ~\vert~ a \leq x \leq b \},
$$
always considered with the usual adjacency relation $c_1$ in $\mathbb{Z}$.

There are two natural choices for the adjacency relations on the Cartesian product of digital images:

\begin{definition}(\cite{boxe17})
Let $(X,\kappa)$ and $(Y,\delta)$ be digital images. Then there are two natural ways to define an adjacency relation on the set $X\times Y$, which we denote $\kappa\times_1 \delta$ and $\kappa \times_2 \delta$.
For pairs $(a,b),(c,d)\in X\times Y$, we define:
\begin{align*}
&(a,b)\sim_{\kappa\times_1 \delta} (c,d) \text{ if and only if: } a\sim_\kappa c \text{ and } b=d, \text{ or } a=c \text{ and } b\sim_\delta d. \\
&(a,b)\sim_{\kappa\times_2 \delta} (c,d) \text{ if and only if: } a\sim_\kappa c \text{ and } b\sim_\delta d.
\end{align*}
In \cite{boxe17} these are called the two \emph{normal product adjacencies}. In the graph theory literature, the graph product formed using $\times_1$ is often called the \emph{box product}, while the product using $\times_2$ is called the \emph{categorical product}.

When we are not indicating the adjacency relations $\kappa$ and $\delta$ explicitly, we will write $X\times_i Y$ to denote $(X\times Y, \kappa \times_i \delta)$.
\end{definition}

The two products $\times_1$ and $\times_2$ correspond naturally to the two adjacencies $c_1$ and $c_2$ on $\Z^2$ as follows: If we take a two-fold product of the usual adjacency for $\Z$, then the induced adjacency on $\Z \times_1 \Z$ is $c_1$, and the adjacency on $\Z\times_2 \Z$ is $c_2$.

These two products also lead immediately to two different homotopy relations:

\begin{definition} (\cite {boxe17, stae21}) \label{homotopydef}
Let $X$ and $Y$ be digital images, and let $f, g : X \to Y$ be digital continuous functions.
For some $i\in \{1,2\}$, an \emph{$i$-homotopy from $f$ to $g$} is a continuous map
\[
H : X \times_i [0,m]_\Z \to Y
\]
such that $H(x,0)=f(x)$ and $H(x,m)=g(x)$ for all $x\in X$.
If such an $H$ exists, we write $f\simeq_i g$, and we say $f$ and $g$ are \emph{$i$-homotopic}.
When $m=1$, we say the homotopy is \emph{single step}.
\end{definition}

Since any adjacency in $X\times_1 [0,m]_\Z$ is also an adjacency in $X\times_2 [0,m]_\Z$, every $2$-homotopy is automatically a $1$-homotopy, and so $f\simeq_2 g$ implies $f\simeq_1 g$.

The following appears as Propositions 1.4 and 2.4 of \cite{stae21}.
\begin{lem}\label{one-step-htp}
Let $X,Y$ be digital images and $f,g:X\to Y$ be continuous. Then:
\begin{itemize}
\item $f\simeq_1 g$ in a single step if and only $f(x)\sim g(x)$ for all $x\in X$.
\item $f\simeq_2 g$ in a single step if and only if $x_1\sim x_2$ implies $f(x_1)\sim g(x_2)$ for all $x_1,x_2 \in X$.
\end{itemize}
\end{lem}

The simplest possible homotopy is one in a single step which moves only one point:
\begin{definition}
For maps $f,g:X\to Y$ and some point $a\in X$, a single-step $i$-homotopy $f\simeq_i g$ is called a \emph{spider move at $a$} when $f(x)=g(x)$ for all $x\neq a$.
\end{definition}

It is easy to show that any spider move is automatically a $2$-homotopy (and thus also a $1$-homotopy).

The term \emph{spider move} was coined in \cite{cs21}, an analogy to a spider which walks by moving only one  leg at a time. In \cite{stae21} this type of homotopy was called a \emph{punctuated homotopy}.
The following was proved independently in \cite{cs21} and \cite{stae21}:
\begin{lem}[The Spider Move Lemma]\label{SML}
For finite digital images $X,Y$ and maps $f,g:X\to Y$, let $f\simeq_2 g$. Then there is a $2$-homotopy $H:X\times_2 [0,k]_\Z \to Y$ from $f$ to $g$ such that, for each induced map $H_t(x) = H(x,t)$, the single-step homotopy $H_t\simeq_2 H_{t+1}$ is a spider move.
\end{lem}

A homotopy $H$ as in Lemma \ref{SML} is referred to as a \emph{homotopy through spider moves}.
%
%

The following lemma gives an easy way to check if $\id_X$ can be changed by a spider move.
In a digital image $(X,\kappa)$ with some point $x\in X$, we define $N(x)$ as the set of all points adjacent to $x$ (this set will always include $x$ itself). When we need to indicate the image or adjacency, we will write this set as $N_X(x)$.
\begin{lem}\label{2htp-nbhd}
Let $f:X\to X$ be continuous. Then there is a map $2$-homotopic to the identity with $f\neq \id_X$ if and only if there is some pair of points $a,b\in X$ with $N(a)\subseteq N(b)$.
\end{lem}
\begin{proof}
First assume there is a map $2$-homotopic to the identity with $f\neq \id_X$.
Let $g:X\to X$ be the first non-identity stage of the $2$-homotopy from $f$ to $\id_X$, and let $a\in X$ be some point with $b=g(a)\neq a$. We will show that $N(a)\subseteq N(b)$. Take any $c\in N(a)$, and since $c\sim a$ and $\id_X\simeq_2 g$ in one step we have $c \sim g(a)=b$ by Lemma \ref{one-step-htp}, and thus $c\in N(b)$ as desired.

For the converse, assume there are points with $N(a)\subseteq N(b)$. Then let:
\[ f(x) = \begin{cases} x &$ if $ x\neq a \\
b &$ if $ x=a. \end{cases} \]
It is routine to check that $f$ is continuous and $f\simeq_2 \id_X$.
\end{proof}


%

As in classical topology, digital images $X$ and $Y$ are called \emph{$i$-homotopy equivalent} when there are continuous maps $f:X\to Y$ and $g:Y\to X$ with $f\circ g \simeq_i \id_Y$ and $g\circ f \simeq_i \id_X$.


\begin{lem}\label{htp-equiv-image}
Let $X$ be a digital image with a map $f:X\to X$ such that $f\simeq_i \id_X$. Then $X$ and $f(X)$ are $i$-homotopy equivalent.
\end{lem}
\begin{proof}
Let $Y := f(X)$ and let $\bar f : X \rightarrow Y$ be the continuous map obtained by restricting the codomain of $f : X \rightarrow X$ to $Y$; that is, $\bar f(x) = f(x)$ for all $x \in X$.
Then we have
$$
i_Y \circ \bar f = f \simeq_i \id_X
$$
by our hypothesis, where $i_Y : Y \hookrightarrow X$ is the inclusion map.

The map $\bar f \circ i_Y: Y \to Y$ is the same as the restriction of $f$ to $Y$. Restricting the $i$-homotopy of $f$ to $\id_X$ shows that $\bar f\circ i_Y \simeq_1 \id_Y$.

Therefore, $\bar f : X \rightarrow Y$ is an $i$-homotopy equivalence whose homotopy-inverse is the inclusion map $i_Y : Y \hookrightarrow X$.
\end{proof}
%

Many of our results will focus on reducing a digital image by an $i$-homotopy equivalence to a proper subset of itself. For 2-homotopy, this can be done if we have certain neighborhood containments.

\begin{lem}\label{htp-equiv-remove-pt}
Let $X$ be a digital image and let $p,q\in X$ be points with $N(p)\subseteq N(q)$. Then $X-\{p\}$ is 2-homotopy equivalent to $X$.
\end{lem}
\begin{proof}
Since $N(p) \subseteq N(q)$, the construction of  Lemma \ref{2htp-nbhd} gives a surjective map $f:X\to X-\{p\}$ such that $f$ is $2$-homotopic to $\id_X$ in a single step. 
\end{proof}

Work in \cite{hmps15} includes a version of Lemma \ref{htp-equiv-remove-pt} which removes an entire path of points from $X$ under the appropriate hypotheses. This result holds only for $1$-homotopy:
\begin{lem}[Theorem 4.6 of \cite{hmps15}]\label{htp-equiv-remove-path}
Let $X$ be a digital image and let $\gamma,\delta:[0,k]_\Z \to X$ be injective paths.
Assume that $N(\gamma(i)) \subset N(\delta(i)) \cup \{\gamma(i-1), \gamma(i+1)\}$ for each $0<i<k$, and that $N(\gamma(0)) \subset N(\delta(0)) \cup \{\gamma(1)\}$ and $N(\gamma(k)) \subset N(\delta(k)) \cup \{\gamma(k-1)\}$. Then $X - \gamma([0,k]_\Z)$ is 1-homotopy equivalent to $X$.
\end{lem}

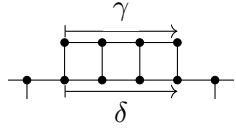
\begin{figure}
\[
\begin{tikzpicture}
\begin{scope}
\clip (-.5,.5) rectangle (5.5,3);
\draw[edge] (-1, 1)--(0, 1);
\draw[edge] (4, 1)--(4, 2)--(3, 2);
\draw[edge] (5, 0)--(5, 1)--(6, 1);
\draw[edge] (3, 1)--(4, 1)--(5, 1);
\draw[edge] (2, 1)--(3, 1)--(3, 2)--(2, 2);
\draw[edge] (0, 0)--(0, 1)--(1, 1);
\draw[edge] (1, 2)--(1, 1)--(2, 1)--(2, 2)--(1, 2);
\foreach \x/\y in {0/0,5/0,-1/1,0/1,1/1,2/1,3/1,4/1,5/1,6/1,1/2,2/2,3/2,4/2} {
  \node[vertex] at (\x,\y) {};
}
\end{scope}
\draw[{Bar}->] (1,2.3) -- (4,2.3) node[above,pos=.5,fill=white] {$\gamma$};
\draw[{Bar}->] (1,.7) -- (4,.7) node[below,pos=.5,fill=white] {$\delta$};
\end{tikzpicture}
\]
\caption{Demonstration of Lemma \ref{htp-equiv-remove-path}. The path $\gamma$ can be removed from $X$ by a 1-homotopy equivalence.\label{remove-path-fig}}
\end{figure}

The hypotheses above are complicated, but they indicate that the points of $\gamma$ and $\delta$ are adjacent in a natural way. See Figure \ref{remove-path-fig} for a typical situation where the Lemma applies.

%
%
%
%

\bigskip

\section{Complements and digital pictures in $\Z^2$}\label{complementsandpicture}
\begin{figure}
\[
\vcbox{
\begin{tikzpicture}
\draw[edge] (1, 0)--(0, 0)--(0, 1)--(1, 1)--(1, 0);
\foreach \x/\y in {0/0,0/1,1/0,1/1} {
  \node[vertex] at (\x,\y) {};
}
\end{tikzpicture}
}
\qquad
\vcbox{
\begin{tikzpicture}
\draw[edge] (-1, 0)--(0, -1)--(1, 0)--(0, 1)--(-1, 0);
\foreach \x/\y in {1/0,0/-1,-1/0,0/1} {
  \node[vertex] at (\x,\y) {};
}
\end{tikzpicture}
}
\]
\caption{The square $(C_4^\square,c_1)$ and the diamond $(C_4^\diamond,c_2)$ are isomorphic as graphs, but the diamond has ``a hole'' while the square does not.\label{C4fig}}
\end{figure}
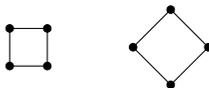

The example of $C_4$, the simple cycle of 4 points, provides an immediate case in which the notion of digital isomorphism does not seem to match our intuitive notion of topological equivalence. The cycle $C_4$ may be embedded in $\Z^2$ in two different ways which we call \emph{the square} and \emph{the diamond} (see Figure \ref{C4fig}):
\[ \begin{split}
(C_4^\square,c_1) &= (\{(0,0),(0,1),(1,0),(1,1) \} ,c_1)\\
(C_4^\diamond,c_2) &= (\{(1,0),(0,-1),(-1,0),(0,1)\},c_2)
\end{split}
\]

These two digital images are isomorphic as graphs, but on an intuitive level they seem to be topologically different because the diamond has a ``hole'' in the middle while the square does not.
Digital $i$-homotopy equivalence offers a more topologically sophisticated notion of equivalence, but since the square and diamond are isomorphic as graphs, they are automatically $i$-homotopy equivalent for both $i\in \{1,2\}$.

\begin{figure}
\[
\vcbox{
\begin{tikzpicture}
\foreach \x in {-1,0,1,2} {
  \foreach \y in {-1,0,1,2} {
    \node[point] at (\x,\y) {};
}}
\foreach \x/\y in {0/0,0/1,1/0,1/1} {
  \node[vertex] at (\x,\y) {};
}
\draw[edge] (1, 1)--(0, 1) (1, 1)--(1, 0) ;
\draw[edge] (1, 0)--(0, 0) ;
\draw[edge] (0, 1)--(0, 0) ;
\begin{scope}
\clip (-1.25, -1.25) rectangle (2.25, 2.25);
\draw[edge] (3, 3)--(2, 2) (3, 3)--(2, 3) (3, 3)--(3, 2) ;
\draw[edge] (3, 2)--(2, 1) (3, 2)--(2, 2) (3, 2)--(2, 3) (3, 2)--(3, 1) ;
\draw[edge] (3, 1)--(2, 0) (3, 1)--(2, 1) (3, 1)--(2, 2) (3, 1)--(3, 0) ;
\draw[edge] (3, 0)--(2, -1) (3, 0)--(2, 0) (3, 0)--(2, 1) (3, 0)--(3, -1) ;
\draw[edge] (3, -1)--(2, -2) (3, -1)--(2, -1) (3, -1)--(2, 0) (3, -1)--(3, -2) ;
\draw[edge] (3, -2)--(2, -2) (3, -2)--(2, -1) ;
\draw[edge] (2, 3)--(1, 2) (2, 3)--(1, 3) (2, 3)--(2, 2) ;
\draw[edge] (2, 2)--(1, 2) (2, 2)--(1, 3) (2, 2)--(2, 1) ;
\draw[edge] (2, 1)--(1, 2) (2, 1)--(2, 0) ;
\draw[edge] (2, 0)--(1, -1) (2, 0)--(2, -1) ;
\draw[edge] (2, -1)--(1, -2) (2, -1)--(1, -1) (2, -1)--(2, -2) ;
\draw[edge] (2, -2)--(1, -2) (2, -2)--(1, -1) ;
\draw[edge] (1, 3)--(0, 2) (1, 3)--(0, 3) (1, 3)--(1, 2) ;
\draw[edge] (1, 2)--(0, 2) (1, 2)--(0, 3) ;
\draw[edge] (1, -1)--(0, -2) (1, -1)--(0, -1) (1, -1)--(1, -2) ;
\draw[edge] (1, -2)--(0, -2) (1, -2)--(0, -1) ;
\draw[edge] (0, 3)--(-1, 2) (0, 3)--(-1, 3) (0, 3)--(0, 2) ;
\draw[edge] (0, 2)--(-1, 1) (0, 2)--(-1, 2) (0, 2)--(-1, 3) ;
\draw[edge] (0, -1)--(-1, -2) (0, -1)--(-1, -1) (0, -1)--(-1, 0) (0, -1)--(0, -2) ;
\draw[edge] (0, -2)--(-1, -2) (0, -2)--(-1, -1) ;
\draw[edge] (-1, 3)--(-2, 2) (-1, 3)--(-2, 3) (-1, 3)--(-1, 2) ;
\draw[edge] (-1, 2)--(-2, 1) (-1, 2)--(-2, 2) (-1, 2)--(-2, 3) (-1, 2)--(-1, 1) ;
\draw[edge] (-1, 1)--(-2, 0) (-1, 1)--(-2, 1) (-1, 1)--(-2, 2) (-1, 1)--(-1, 0) ;
\draw[edge] (-1, 0)--(-2, -1) (-1, 0)--(-2, 0) (-1, 0)--(-2, 1) (-1, 0)--(-1, -1) ;
\draw[edge] (-1, -1)--(-2, -2) (-1, -1)--(-2, -1) (-1, -1)--(-2, 0) (-1, -1)--(-1, -2) ;
\draw[edge] (-1, -2)--(-2, -2) (-1, -2)--(-2, -1) ;
\draw[edge] (-2, 3)--(-2, 2) ;
\draw[edge] (-2, 2)--(-2, 1) ;
\draw[edge] (-2, 1)--(-2, 0) ;
\draw[edge] (-2, 0)--(-2, -1) ;
\draw[edge] (-2, -1)--(-2, -2) ;
\foreach \x/\y in {-2/-2,-2/-1,-2/0,-2/1,-2/2,-2/3,-1/-2,-1/-1,-1/0,-1/1,-1/2,-1/3,0/-2,0/-1,0/2,0/3,1/-2,1/-1,1/2,1/3,2/-2,2/-1,2/0,2/1,2/2,2/3,3/-2,3/-1,3/0,3/1,3/2,3/3} {
  \node[comp-vertex] at (\x,\y) {};
}
\end{scope}
\end{tikzpicture}
}
\qquad
\vcbox{
\begin{tikzpicture}
\foreach \x in {-2,-1,0,1,2} {
  \foreach \y in {-2,-1,0,1,2} {
    \node[point] at (\x,\y) {};
}}
\foreach \x/\y in {1/0,0/-1,-1/0,0/1} {
  \node[vertex] at (\x,\y) {};
}
\draw[edge] (0, 1)--(1, 0) (0, 1)--(-1, 0) ;
\draw[edge] (-1, 0)--(0, -1) ;
\draw[edge] (0, -1)--(1, 0) ;
\begin{scope}
\clip (-2.25, -2.25) rectangle (2.25, 2.25);
\draw[edge] (3, 3)--(2, 3) (3, 3)--(3, 2) ;
\draw[edge] (3, 2)--(2, 2) (3, 2)--(3, 1) ;
\draw[edge] (3, 1)--(2, 1) (3, 1)--(3, 0) ;
\draw[edge] (3, 0)--(2, 0) (3, 0)--(3, -1) ;
\draw[edge] (3, -1)--(2, -1) (3, -1)--(3, -2) ;
\draw[edge] (3, -2)--(2, -2) (3, -2)--(3, -3) ;
\draw[edge] (3, -3)--(2, -3) ;
\draw[edge] (2, 3)--(1, 3) (2, 3)--(2, 2) ;
\draw[edge] (2, 2)--(1, 2) (2, 2)--(2, 1) ;
\draw[edge] (2, 1)--(1, 1) (2, 1)--(2, 0) ;
\draw[edge] (2, 0)--(2, -1) ;
\draw[edge] (2, -1)--(1, -1) (2, -1)--(2, -2) ;
\draw[edge] (2, -2)--(1, -2) (2, -2)--(2, -3) ;
\draw[edge] (2, -3)--(1, -3) ;
\draw[edge] (1, 3)--(0, 3) (1, 3)--(1, 2) ;
\draw[edge] (1, 2)--(0, 2) (1, 2)--(1, 1) ;
\draw[edge] (1, -1)--(1, -2) ;
\draw[edge] (1, -2)--(0, -2) (1, -2)--(1, -3) ;
\draw[edge] (1, -3)--(0, -3) ;
\draw[edge] (0, 3)--(-1, 3) (0, 3)--(0, 2) ;
\draw[edge] (0, 2)--(-1, 2) ;
\draw[edge] (0, -2)--(-1, -2) (0, -2)--(0, -3) ;
\draw[edge] (0, -3)--(-1, -3) ;
\draw[edge] (-1, 3)--(-2, 3) (-1, 3)--(-1, 2) ;
\draw[edge] (-1, 2)--(-2, 2) (-1, 2)--(-1, 1) ;
\draw[edge] (-1, 1)--(-2, 1) ;
\draw[edge] (-1, -1)--(-2, -1) (-1, -1)--(-1, -2) ;
\draw[edge] (-1, -2)--(-2, -2) (-1, -2)--(-1, -3) ;
\draw[edge] (-1, -3)--(-2, -3) ;
\draw[edge] (-2, 3)--(-3, 3) (-2, 3)--(-2, 2) ;
\draw[edge] (-2, 2)--(-3, 2) (-2, 2)--(-2, 1) ;
\draw[edge] (-2, 1)--(-3, 1) (-2, 1)--(-2, 0) ;
\draw[edge] (-2, 0)--(-3, 0) (-2, 0)--(-2, -1) ;
\draw[edge] (-2, -1)--(-3, -1) (-2, -1)--(-2, -2) ;
\draw[edge] (-2, -2)--(-3, -2) (-2, -2)--(-2, -3) ;
\draw[edge] (-2, -3)--(-3, -3) ;
\draw[edge] (-3, 3)--(-3, 2) ;
\draw[edge] (-3, 2)--(-3, 1) ;
\draw[edge] (-3, 1)--(-3, 0) ;
\draw[edge] (-3, 0)--(-3, -1) ;
\draw[edge] (-3, -1)--(-3, -2) ;
\draw[edge] (-3, -2)--(-3, -3) ;
\foreach \x/\y in {-3/-3,-3/-2,-3/-1,-3/0,-3/1,-3/2,-3/3,-2/-3,-2/-2,-2/-1,-2/0,-2/1,-2/2,-2/3,-1/-3,-1/-2,-1/-1,-1/1,-1/2,-1/3,0/-3,0/-2,0/0,0/2,0/3,1/-3,1/-2,1/-1,1/1,1/2,1/3,2/-3,2/-2,2/-1,2/0,2/1,2/2,2/3,3/-3,3/-2,3/-1,3/0,3/1,3/2,3/3} {
  \node[comp-vertex] at (\x,\y) {};
}
\end{scope}
\end{tikzpicture}
}
\]
\caption{The square $(C_4^\square,c_1,c_2)$ and the diamond $(C_4^\diamond,c_2,c_1)$, as digital pictures, are not isomorphic.\label{C4complementsfig}}
\end{figure}
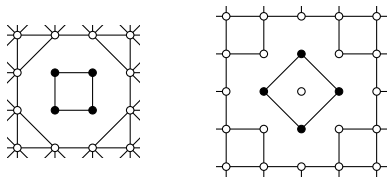

The situation becomes clearer when we consider the square and the diamond each together with their complements, using the complementary adjacencies. This is shown in Figure \ref{C4complementsfig}. The points of the digital image are shown in the black vertices, and the ponts of the complement are the white vertices.

We see immediately that the complements of the square and the diamond are not isomorphic, since the square has connected complement while the diamond does not. From this point of view, the topological properties of a digital image are best expressed when we simultaneously consider its complement as a digital image.

Let $\bar X$ denote the set-theoretic complement of a digital image $X$; that is, $\bar X := X^c \subsetneq \mathbb{Z}^2$.

Following terminology used in the classical digital topology literature (see \cite{kr89}), a digital image with its complement is called a \emph{digital picture}. One of the fundamental characteristics of the work in digital pictures is the decision to use different adjacency relations for $X$ and for its complement $\bar X$, as we did in Figure \ref{C4complementsfig}. This paper deals exclusively with digital pictures in $\Z^2$, so we make the following definition:
\begin{definition}
A \emph{digital picture} is a triple $(X,c_i,c_j)$ where $X\subsetneq \Z^2$ is a finite set and $\{i,j\} = \{1,2\}$. In this context, the set $X$ is called the set of \emph{black pixels}, and is viewed as a digital image with adjacency $c_i$. The complement $\bar X = \Z^2-X$ is called the set of \emph{white pixels}, and it is viewed as a digital image with adjacency $c_j$.
\end{definition}

The use of two different adjacencies for $X$ and $\bar X$ is counterintuitive but useful: it is the proper way to recognize the difference between the square and the diamond. It also gives the correct setting to formulate the digital Jordan Curve Theorem:
\begin{thm}[Rosenfeld 1979 \cite{rose79}]
Let $(X,c_i,c_j)$ be a digital picture where $X\subset \Z^2$ is a digital simple closed $c_i$-curve. Assume that $X$ has at least 4 points when $i=1$, and at least 8 points when $i=2$. Then $(\bar X,c_j)$ has two $c_j$-components: one finite, the other infinite.
\end{thm}

Above, a ``digital simple closed $c_i$-curve'' is a $c_i$-continuous map $\gamma: [0,k]_\Z \to X$ with $\gamma(i)\not \sim_{c_i} \gamma(j)$ for all $i,j$ except $\gamma(i) \sim_{c_i} \gamma(i+1)$ and $\gamma(0)=\gamma(k)$.



Our main goal for the paper is to explore the notion of homotopy equivalence for digital pictures. We begin with a simpler isomorphism relation.

\begin{definition}
Let $(X,\kappa,\bar \kappa)$ and $(Y,\delta,\bar \delta)$ be digital pictures. We say they are \emph{isomorphic} when the digital images $(X,\kappa)$ and $(Y,\delta)$ are isomorphic, and $(\bar X,\bar \kappa)$ and $(\bar Y,\bar \delta)$ are isomorphic.
\end{definition}

Thus $(C_4^\square,c_1)$ and $(C_4^\diamond,c_2)$ are isomorphic as digital images, but $(C_4^\square,c_1,c_2)$ and $(C_4^\diamond,c_2,c_1)$ are not isomorphic as digital pictures. In fact no digital pictures $(X,c_1,c_2)$ and $(Y,c_2,c_1)$ can be isomorphic because of the different adjacencies used in the complement. In this case $(\bar X,c_2)$ is an infinite graph with all but finitely many vertices having degree 8, while $(\bar Y, c_1)$ is an infinite graph with all but finitely many vertices having degree 4, so they are not isomorphic.

We will define a homotopy equivalence relation for digital pictures in a similar way to isomorphism, requiring that both the sets and their complements be homotopy equivalent. There is a choice to consider, however, regarding the type of homotopy to use: either 1-homotopy or 2-homotopy.

Because of the Spider Move Lemma, the 2-homotopy relation is very restrictive in a digital image $(X,c_1)$. For example the condition $N(a)\subseteq N(b)$ used in Lemma \ref{2htp-nbhd} can only hold in a $c_1$-digital image when $a$ is adjacent to $b$ and no other points. In practice this makes 2-homotopy an overly restrictive relation when using $c_1$-adjacency.

Unsurprisingly, the work in \cite{los21} and following papers by the same authors, which focus exclusively on $2$-homotopy, never uses $c_1$-adjacency. Those papers achieving interesting results by always using $2$-homotopy, and $c_n$-adjacency in $\Z^n$.

The field of A-Theory, which exclusively uses $1$-homotopy, has had great success with ``cubical'' constructions, in particular discrete cubical homology theory. In this case we see that $c_1$-adjacency (which naturally makes a digital image into a cubical complex) gives strong results for constructions using 1-homotopy.

Thus for digital images in $\Z^2$, we arrive at the point of view that $1$-homotopy is the natural choice when working with $c_1$-adjacency, and $2$-homotopy is the natural choice when working with $c_2$-adjacency. That is, when considering a digital image $(X,c_i)$ for $i\in \{1,2\}$, we will use $i$-homotopy when considering maps on $X$. For a digital picture $(X,c_i,c_j)$, we use $i$-homotopy in $X$ and $j$-homotopy in $\bar X$. As we will see, the main results in our later sections require this use of two different types of homotopy in $X$ and $\bar X$.

\begin{definition}
Let $(X,c_i,c_j)$ and $(Y,c_i,c_j)$ be digital pictures. We say they are \emph{homotopy equivalent}, or they have the same \emph{homotopy type}, when the digital images $(X,c_i)$ and $(Y,c_i)$ are $i$-homotopy equivalent, and $(\bar X,c_j)$ and $(\bar Y,c_j)$ are $j$-homotopy equivalent.
\end{definition}

Note that we do not define homotopy equivalence of digital pictures $(X,c_i,c_j)$ and $(Y,c_j,c_i)$. This would require more choices about which type of homotopy equivalences to use for the white and black points.

\bigskip

\section{The outer perimeter, and digital rectangles}\label{digitalrectangles}
This paper aims to explore the homotopy equivalence relation for digital pictures.
As a simple starting point, we discuss homotopy equivalences between digital pictures consisting of rectangles in $\Z^2$. As digital images, these are all $i$-contractible for both $i\in \{1,2\}$, and thus they all have the same $i$-homotopy type. But as digital pictures we will see that rectangles of different sizes will generally have different homotopy type.

To help distinguish the homotopy type of rectangles of different sizes, we will adapt the loop-counting invariant discussed in \cite{hmps15}, which will give a simple numerical homotopy-type invariant for digital pictures.
\begin{definition}
For a digital picture $(X,c_i,c_j)$, let $O(X,c_i,c_j)$ be the smallest possible length of a non-$j$-contractible loop in the infinite component of $\bar X$. We call this the number the \emph{outer perimeter} of $(X,c_i,c_j)$. If no non-$j$-contractible loop exists in the infinite component of $\bar X$, then we define $O(X,c_i,c_j)=1$.
\end{definition}

Next we show that the outer perimeter is a homotopy invariant for digital pictures.

\begin{prop}
For any digital picture $(X,c_i,c_j)$, the outer perimeter $O(X,c_i,c_j)$ is a positive integer, and $O(X,c_i,c_j)=1$ only when $X$ is empty.
\end{prop}
\begin{proof}
Clearly $O(X,c_i,c_j)$ is a positive integer because it is a length of some loop in $\bar X$. For the second statement, observe that the smallest non-1-contractible $c_1$-loop which surrounds a nonempty set has length 8, while the smallest non-2-contractible $c_2$-loop which surrounds a nonempty set has length 4. Thus in fact we always have $O(X,c_i,c_j) \ge 4$ when $X$ is nonempty.
\end{proof}

\begin{prop}
The outer perimeter is a homotopy-type invariant. That is, if $(X,c_i,c_j)$ and $(Y,c_i,c_j)$ are homotopy equivalent digital pictures, then $O(X,c_i,c_j) = O(Y,c_i,c_j)$.
\end{prop}
\begin{proof}
The invariant $O(X,c_i,c_j)$ is a simplified form of the family of loop counting invariants $L_m(X,\kappa)$ defined in \cite{hmps15}. This $L_m(X,\kappa)$ is defined as the number of non-contractible loops in $X$ of length $m$, counted modulo a natural homotopy relation. Our invariant $O(X,c_i,c_j)$ is simply the smallest $m$ for which $L_m(\bar X,c_j)$ is nonzero. Work in \cite{hmps15} shows that all of the numbers $L_m(X,\kappa)$ for various $m$ are $i$-homotopy type invariants. (The proof in \cite{hmps15} is written for $1$-homotopy, but the same argument holds exactly as written for 2-homotopy.)

Thus if $(X,c_i,c_j)$ and $(Y,c_i,c_j)$ are homotopy equivalent, then their entire sequences of numbers $(L_m(\bar X,c_j))_{m=1}^\infty$ and $(L_m(\bar Y,c_j))_{m=1}^\infty$ are equal, and so we will have $O(X,c_i,c_j) = O(Y,c_i,c_j)$.
\end{proof}

Let $I_{n,m}$ be the digital rectangle $I_{n,m} = [1,n]_\Z \times [1,m]_\Z$, and we wish to describe the homotopy type of $(I_{n,m},c_i,c_j)$.

First we discuss the isomorphism type of rectangles:
\begin{thm}\label{rectangle-iso}
Let $\{i,j\} = \{1,2\}$. Two digital rectangles $(I_{n,m},c_i,c_j)$ and $(I_{k,l},c_i,c_j)$ are isomorphic if and only if $\{n,m\} = \{k,l\}$.
\end{thm}
\begin{proof}
Clearly if $\{n,m\} = \{k,l\}$, then $(I_{n,m},c_1,c_2)$ and $(I_{k,l},c_1,c_2)$ are isomorphic either by an identity or reflection map. For the converse, assume that $(I_{n,m},c_i,c_j)$ and $(I_{k,l},c_i,c_j)$ are isomorphic. This is now a purely graph theoretic argument.

For $i=1$, the graph $I_{n,m}$ has $nm$ total vertices, of which $2(n-2)+2(m-2) = 2n+2m-8$ have degree 3. Thus if $I_{n,m}$ and $I_{k,l}$ are isomorphic then we must have $nm=kl$ and $n+m=k+l$. Thus we have an equality of quadratic polynomials:
\[ x^2 - (n+m)x + nm = x^2 - (k+l)x+kl, \]
and factoring gives $(x-n)(x-m) = (x-k)(x-l)$ and thus $\{n,m\} = \{k,l\}$.

For $i=2$ we argue similarly, this time counting total vertices, and vertices of degree 5, of which there are again $2n+2m-8$. We again conclude that $\{n,m\} = \{k,l\}$.
\end{proof}

For homotopy equivalence rather than isomorphism, we have a weaker result.

\begin{thm}\label{rectangle-perimeter}
Let $\{i,j\} = \{1,2\}$. If the digital pictures $(I_{n,m},c_i,c_j)$ and $(I_{k,l},c_i, c_j)$ are homotopy equivalent, then $n+m= k+l$. That is, the rectangles have the same outer perimeter.
\end{thm}
\begin{proof}

Assuming that $(I_{n,m},c_i,c_j)$ and $(I_{k,l},c_i,c_j)$ are homotopy equivalent, we must have $O(I_{n,m},c_i,c_j) = O(I_{k,l},c_i,c_j)$. It is easy to verify that $O(I_{n,m},c_1,c_2) = 2n+2m$, and $O(I_{n,m},c_2,c_1) = 2n+2m+4$. Thus $O(I_{n,m},c_i,c_j) = O(I_{k,l},c_i,c_j)$ will imply that $n+m=k+l$ as desired.
\end{proof}

Let $I_n = [1,n]_\Z \times \{1\}$ be a digital interval in $\Z^2$. Using $m=l=1$ in the theorem above gives:
\begin{thm}
Let $\{i,j\} = \{1,2\}$. The digital pictures $(I_n,c_i,c_j)$ and $(I_k,c_i,c_j)$ are homotopy equivalent if and only if $n=k$.
\end{thm}

The conclusion of Theorem \ref{rectangle-perimeter} says that homotopy equivalent rectangles must have the same perimeter. We do not know if this statement can be strengthened.

\begin{quest}
If $(I_{n,m},c_i,c_j)$ and $(I_{k,l},c_i,c_j)$ are homotopy equivalent, then must they be isomorphic?
\end{quest}
As a specific example, $(I_{3,1},c_i,c_j)$ and $(I_{2,2},c_i, c_j)$ are not isomorphic because $I_{3,1}$ and $I_{2,2}$ have different cardinalities. But $3+1=2+2$, so they do not violate the conclusion of Theorem \ref{rectangle-perimeter}. We do not know if these are homotopy equivalent or not.

For non-rectangles, it is possible for $(X,c_i, c_j)$ and $(Y,c_i, c_j)$ to be homotopy equivalent but not isomorphic.

\begin{exa}\label{bowtie-ex}
We consider the digital pictures $(X,c_1,c_2)$ and $(I_{4,2},c_1,c_2)$ where
\[ X = \{(1,2),(1,1),(2,1),(3,1),(4,1),(4,2) \}, \]
shown in Figure \ref{bowtiefig}. Then $X$ and $I_{4,2}$ have different cardinalities, so $(X,c_1,c_2)$ and $(I_{4,2},c_1,c_2)$ are not isomorphic. But we argue that $(X,c_1,c_2)$ and $(I_{4,2},c_1,c_2)$ are homotopy equivalent.
\begin{figure}
\[
\begin{tikzpicture}
\foreach \x/\y in {0/0,1/0,2/0,3/0,0/1,3/1} {
  \node[vertex] at (\x,\y) {};
}
\draw[edge] (3, 1)--(3, 0) ;
\draw[edge] (0, 1)--(0, 0) ;
\draw[edge] (3, 0)--(2, 0) ;
\draw[edge] (2, 0)--(1, 0) ;
\draw[edge] (1, 0)--(0, 0) ;
\begin{scope}
\clip (-1.25, -1.25) rectangle (4.25, 2.25);
\draw[edge] (5, 3)--(4, 2) (5, 3)--(4, 3) (5, 3)--(5, 2) ;
\draw[edge] (5, 2)--(4, 1) (5, 2)--(4, 2) (5, 2)--(4, 3) (5, 2)--(5, 1) ;
\draw[edge] (5, 1)--(4, 0) (5, 1)--(4, 1) (5, 1)--(4, 2) (5, 1)--(5, 0) ;
\draw[edge] (5, 0)--(4, -1) (5, 0)--(4, 0) (5, 0)--(4, 1) (5, 0)--(5, -1) ;
\draw[edge] (5, -1)--(4, -2) (5, -1)--(4, -1) (5, -1)--(4, 0) (5, -1)--(5, -2) ;
\draw[edge] (5, -2)--(4, -2) (5, -2)--(4, -1) ;
\draw[edge] (4, 3)--(3, 2) (4, 3)--(3, 3) (4, 3)--(4, 2) ;
\draw[edge] (4, 2)--(3, 2) (4, 2)--(3, 3) (4, 2)--(4, 1) ;
\draw[edge] (4, 1)--(3, 2) (4, 1)--(4, 0) ;
\draw[edge] (4, 0)--(3, -1) (4, 0)--(4, -1) ;
\draw[edge] (4, -1)--(3, -2) (4, -1)--(3, -1) (4, -1)--(4, -2) ;
\draw[edge] (4, -2)--(3, -2) (4, -2)--(3, -1) ;
\draw[edge] (3, 3)--(2, 2) (3, 3)--(2, 3) (3, 3)--(3, 2) ;
\draw[edge] (3, 2)--(2, 1) (3, 2)--(2, 2) (3, 2)--(2, 3) ;
\draw[edge] (3, -1)--(2, -2) (3, -1)--(2, -1) (3, -1)--(3, -2) ;
\draw[edge] (3, -2)--(2, -2) (3, -2)--(2, -1) ;
\draw[edge] (2, 3)--(1, 2) (2, 3)--(1, 3) (2, 3)--(2, 2) ;
\draw[edge] (2, 2)--(1, 1) (2, 2)--(1, 2) (2, 2)--(1, 3) (2, 2)--(2, 1) ;
\draw[edge] (2, 1)--(1, 1) (2, 1)--(1, 2) ;
\draw[edge] (2, -1)--(1, -2) (2, -1)--(1, -1) (2, -1)--(2, -2) ;
\draw[edge] (2, -2)--(1, -2) (2, -2)--(1, -1) ;
\draw[edge] (1, 3)--(0, 2) (1, 3)--(0, 3) (1, 3)--(1, 2) ;
\draw[edge] (1, 2)--(0, 2) (1, 2)--(0, 3) (1, 2)--(1, 1) ;
\draw[edge] (1, 1)--(0, 2) ;
\draw[edge] (1, -1)--(0, -2) (1, -1)--(0, -1) (1, -1)--(1, -2) ;
\draw[edge] (1, -2)--(0, -2) (1, -2)--(0, -1) ;
\draw[edge] (0, 3)--(-1, 2) (0, 3)--(-1, 3) (0, 3)--(0, 2) ;
\draw[edge] (0, 2)--(-1, 1) (0, 2)--(-1, 2) (0, 2)--(-1, 3) ;
\draw[edge] (0, -1)--(-1, -2) (0, -1)--(-1, -1) (0, -1)--(-1, 0) (0, -1)--(0, -2) ;
\draw[edge] (0, -2)--(-1, -2) (0, -2)--(-1, -1) ;
\draw[edge] (-1, 3)--(-2, 2) (-1, 3)--(-2, 3) (-1, 3)--(-1, 2) ;
\draw[edge] (-1, 2)--(-2, 1) (-1, 2)--(-2, 2) (-1, 2)--(-2, 3) (-1, 2)--(-1, 1) ;
\draw[edge] (-1, 1)--(-2, 0) (-1, 1)--(-2, 1) (-1, 1)--(-2, 2) (-1, 1)--(-1, 0) ;
\draw[edge] (-1, 0)--(-2, -1) (-1, 0)--(-2, 0) (-1, 0)--(-2, 1) (-1, 0)--(-1, -1) ;
\draw[edge] (-1, -1)--(-2, -2) (-1, -1)--(-2, -1) (-1, -1)--(-2, 0) (-1, -1)--(-1, -2) ;
\draw[edge] (-1, -2)--(-2, -2) (-1, -2)--(-2, -1) ;
\draw[edge] (-2, 3)--(-2, 2) ;
\draw[edge] (-2, 2)--(-2, 1) ;
\draw[edge] (-2, 1)--(-2, 0) ;
\draw[edge] (-2, 0)--(-2, -1) ;
\draw[edge] (-2, -1)--(-2, -2) ;
\foreach \x/\y in {-2/-2,-2/-1,-2/0,-2/1,-2/2,-2/3,-1/-2,-1/-1,-1/0,-1/1,-1/2,-1/3,0/-2,0/-1,0/2,0/3,1/-2,1/-1,1/1,1/2,1/3,2/-2,2/-1,2/1,2/2,2/3,3/-2,3/-1,3/2,3/3,4/-2,4/-1,4/0,4/1,4/2,4/3,5/-2,5/-1,5/0,5/1,5/2,5/3} {
  \node[comp-vertex] at (\x,\y) {};
}
\end{scope}
\end{tikzpicture}
\qquad
\begin{tikzpicture}
\foreach \x/\y in {0/0,0/1,1/0,1/1,2/0,2/1,3/0,3/1} {
  \node[vertex] at (\x,\y) {};
}
\draw[edge] (3, 1)--(2, 1) (3, 1)--(3, 0) ;
\draw[edge] (3, 0)--(2, 0) ;
\draw[edge] (2, 1)--(1, 1) (2, 1)--(2, 0) ;
\draw[edge] (2, 0)--(1, 0) ;
\draw[edge] (1, 1)--(0, 1) (1, 1)--(1, 0) ;
\draw[edge] (1, 0)--(0, 0) ;
\draw[edge] (0, 1)--(0, 0) ;
\begin{scope}
\clip (-1.25, -1.25) rectangle (4.25, 2.25);
\draw[edge] (5, 3)--(4, 2) (5, 3)--(4, 3) (5, 3)--(5, 2) ;
\draw[edge] (5, 2)--(4, 1) (5, 2)--(4, 2) (5, 2)--(4, 3) (5, 2)--(5, 1) ;
\draw[edge] (5, 1)--(4, 0) (5, 1)--(4, 1) (5, 1)--(4, 2) (5, 1)--(5, 0) ;
\draw[edge] (5, 0)--(4, -1) (5, 0)--(4, 0) (5, 0)--(4, 1) (5, 0)--(5, -1) ;
\draw[edge] (5, -1)--(4, -2) (5, -1)--(4, -1) (5, -1)--(4, 0) (5, -1)--(5, -2) ;
\draw[edge] (5, -2)--(4, -2) (5, -2)--(4, -1) ;
\draw[edge] (4, 3)--(3, 2) (4, 3)--(3, 3) (4, 3)--(4, 2) ;
\draw[edge] (4, 2)--(3, 2) (4, 2)--(3, 3) (4, 2)--(4, 1) ;
\draw[edge] (4, 1)--(3, 2) (4, 1)--(4, 0) ;
\draw[edge] (4, 0)--(3, -1) (4, 0)--(4, -1) ;
\draw[edge] (4, -1)--(3, -2) (4, -1)--(3, -1) (4, -1)--(4, -2) ;
\draw[edge] (4, -2)--(3, -2) (4, -2)--(3, -1) ;
\draw[edge] (3, 3)--(2, 2) (3, 3)--(2, 3) (3, 3)--(3, 2) ;
\draw[edge] (3, 2)--(2, 2) (3, 2)--(2, 3) ;
\draw[edge] (3, -1)--(2, -2) (3, -1)--(2, -1) (3, -1)--(3, -2) ;
\draw[edge] (3, -2)--(2, -2) (3, -2)--(2, -1) ;
\draw[edge] (2, 3)--(1, 2) (2, 3)--(1, 3) (2, 3)--(2, 2) ;
\draw[edge] (2, 2)--(1, 2) (2, 2)--(1, 3) ;
\draw[edge] (2, -1)--(1, -2) (2, -1)--(1, -1) (2, -1)--(2, -2) ;
\draw[edge] (2, -2)--(1, -2) (2, -2)--(1, -1) ;
\draw[edge] (1, 3)--(0, 2) (1, 3)--(0, 3) (1, 3)--(1, 2) ;
\draw[edge] (1, 2)--(0, 2) (1, 2)--(0, 3) ;
\draw[edge] (1, -1)--(0, -2) (1, -1)--(0, -1) (1, -1)--(1, -2) ;
\draw[edge] (1, -2)--(0, -2) (1, -2)--(0, -1) ;
\draw[edge] (0, 3)--(-1, 2) (0, 3)--(-1, 3) (0, 3)--(0, 2) ;
\draw[edge] (0, 2)--(-1, 1) (0, 2)--(-1, 2) (0, 2)--(-1, 3) ;
\draw[edge] (0, -1)--(-1, -2) (0, -1)--(-1, -1) (0, -1)--(-1, 0) (0, -1)--(0, -2) ;
\draw[edge] (0, -2)--(-1, -2) (0, -2)--(-1, -1) ;
\draw[edge] (-1, 3)--(-2, 2) (-1, 3)--(-2, 3) (-1, 3)--(-1, 2) ;
\draw[edge] (-1, 2)--(-2, 1) (-1, 2)--(-2, 2) (-1, 2)--(-2, 3) (-1, 2)--(-1, 1) ;
\draw[edge] (-1, 1)--(-2, 0) (-1, 1)--(-2, 1) (-1, 1)--(-2, 2) (-1, 1)--(-1, 0) ;
\draw[edge] (-1, 0)--(-2, -1) (-1, 0)--(-2, 0) (-1, 0)--(-2, 1) (-1, 0)--(-1, -1) ;
\draw[edge] (-1, -1)--(-2, -2) (-1, -1)--(-2, -1) (-1, -1)--(-2, 0) (-1, -1)--(-1, -2) ;
\draw[edge] (-1, -2)--(-2, -2) (-1, -2)--(-2, -1) ;
\draw[edge] (-2, 3)--(-2, 2) ;
\draw[edge] (-2, 2)--(-2, 1) ;
\draw[edge] (-2, 1)--(-2, 0) ;
\draw[edge] (-2, 0)--(-2, -1) ;
\draw[edge] (-2, -1)--(-2, -2) ;
\foreach \x/\y in {-2/-2,-2/-1,-2/0,-2/1,-2/2,-2/3,-1/-2,-1/-1,-1/0,-1/1,-1/2,-1/3,0/-2,0/-1,0/2,0/3,1/-2,1/-1,1/2,1/3,2/-2,2/-1,2/2,2/3,3/-2,3/-1,3/2,3/3,4/-2,4/-1,4/0,4/1,4/2,4/3,5/-2,5/-1,5/0,5/1,5/2,5/3} {
  \node[comp-vertex] at (\x,\y) {};
}
\end{scope}
\end{tikzpicture}
\]
\caption{Digital images $(X,c_1,c_2)$ and $(I_{4,2},c_1,c_2)$. These are homotopy equivalent but not isomorphic.\label{bowtiefig}}
\end{figure}
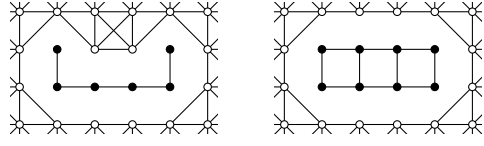

Both $(X,c_1)$ and $(I_{4,2},c_1)$ are $1$-contractible, so they are $1$-homotopy equivalent. For the complements, we have $N_{\bar X}(2,2) \subseteq N_{\bar X}(2,3)$ and $N_{\bar X}(3,2) \subseteq N_{\bar X}(3,3)$. Applying Lemma \ref{htp-equiv-remove-pt} twice shows that $\bar X$ is $2$-homotopy equivalent to $\bar X - \{(2,2),(3,2)\} = \overline{I_{4,2}}$.
\end{exa}

The example above suggests that we can generally use Lemma \ref{htp-equiv-remove-pt} to ``fill in gaps'' in the rows of a digital picture without changing the homotopy type. We explore this in general in the next sections.

\bigskip

\section{The row component graph of a digital image in $\Z^2$}\label{rowcomponent}
We say two sets $A,B\subset \Z^2$ are $c_i$-adjacent when some point of $A$ is $c_i$-adjacent to some point of $B$.

\begin{figure}
\[
\vcbox{
\begin{tikzpicture}[scale=1]
\foreach \x/\y in {0/0,7/0,8/0,9/0,0/1,1/1,2/1,5/1,7/1,9/1,2/2,3/2,4/2,5/2,6/2,7/2,8/2,9/2,5/3,8/3,8/4} {
  \draw[fill=lightgray] (\x,\y) rectangle (\x+.9,\y+.9) {};
}
\end{tikzpicture}
}
\qquad
\vcbox{
\begin{tikzpicture}[scale=1]
\foreach \x/\y in {0/0,2.5/0,0/1,1/1,2/1,3/1,1.5/2,1/3,2/3,2/4} {
 \node[vertex] at (\x,\y) {};
}
\draw (0,0) -- (0,1) -- (1.5,2) -- (1,3);
\draw (1,1) -- (1.5,2) -- (2,3) -- (2,4);
\draw (2,1) -- (1.5,2) -- (3,1) -- (2.5,0) -- (2,1);
\end{tikzpicture}
}
\]
\caption{A $c_1$-digital image and its row component graph.\label{rowgraphfig}}
\end{figure}
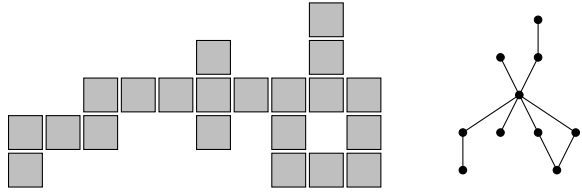

Given a digital image $(X,c_i)$ in $\Z^2$, we define a graph called the \emph{row component graph}. A \emph{row} of $X$ is a set $X\cap (\Z\times \{k\}) \subset X$. The row component graph $R(X,c_i)$ has a vertex for each component of each row of $X$, and two such row components are connected by an edge if they are $c_i$-adjacent. See Figure \ref{rowgraphfig} for a digital image and its row component graph.

From a graph-theoretic point of view, the row component graph $R(X,c_i)$ is obtained by viewing $(X,c_i)$ as a graph and taking the quotient which identifies any two points in the same row component.

\begin{thm}\label{rowgraphconnected}
Let $(X,c_i)$ be a digital picture. The number of $c_i$-components of $X$ equals the number of components of $R(X,c_i)$. In particular, if $(X,c_i)$ is $c_i$-connected, then $R(X,c_i)$ is connected.
\end{thm}
\begin{proof}
This follows immediately from the description of $R(X,c_i)$ as a quotient of $(X,c_i)$. Since this quotient can be achieved by edge contractions, it will preserve the number of components.
\end{proof}

\begin{definition}
For a digital picture $(X,c_i,c_j)$, a finite $c_j$-component of $\bar X$ is called a \emph{hole}.
\end{definition}

\begin{thm}\label{rowgraphacyclic}
Let $(X,c_i,c_j)$ be a digital picture in $\Z^2$. If $X$ has no hole, then $R(X,c_i)$ is an acyclic graph.
\end{thm}
\begin{proof}
We will prove the contrapositive: assume that the graph $R(X,c_i)$ contains a cycle (a loop in the graph with no repeated edges or vertices), and we will show that $X$ has a hole. Let $r_1,\dots,r_k$ be the distinct row components forming a cycle $C$ in $R(X,c_i)$ where $r_i \sim_{c_i} r_{i+1}$ for $i<k$ and $r_1 \sim_{c_i} r_k$.

By taking horizontal paths along each row $r_i$ together with the adjacencies connecting each $r_i$ with $r_{i+1}$, we may use this cycle of row components to build a simple closed $c_i$-curve in $X$ based at a point of $r_1$. Call this curve $\gamma$, and by the digital Jordan Curve Theorem it defines disjoint and non-$c_j$-adjacent ``inside'' and ``outside'' components of $\Z^2$ which we call $A$ and $B$, respectively.

Let $r_l$ be some row component in the cycle $C$ of minimal $y$-coordinate, say $r_l = [a,b] \times \{k\}$. Since $C$ is a cycle, this means that $r_{l-1}$ and $r_{l+1}$ both have $y$-coordinate $k+1$. Since $r_{l-1}$ and $r_{l+1}$ are non-adjacent components of their row, there is some point $(c,k+1)\in \bar X$ with $a\le c \le b$. By this construction, the point $(c,k+1)$ is an element of $A$, the ``inside'' of $\gamma$.

Recall that points in $A$ cannot be connected by a $c_j$-path in $\bar X$ to any point of $B$. Since the infinite component of $\bar X$ is a subset of $B$, this means that $(c,k+1)$ is not in the infinite component of $\bar X$. Thus $\bar X$ has at least one component other than the infinite component, and so $X$ has a hole.
\end{proof}

The Theorem \ref{rowgraphacyclic} above means that the row component graph is tree-like, and the next lemma will show that the leaves of this tree can be removed by homotopy equivalence.

\begin{lem}\label{htp-equiv-remove-row-component}
Let $(X,c_i)$ be a digital image, and $R(X,c_i)$ be its row component graph. If $r\in R(X,c_i)$ is a finite row component which is $c_i$-adjacent to only one other row component, then $(X,c_i)$ is $i$-homotopy equivalent to $(X-r,c_i)$.
\end{lem}
\begin{proof}
Let $s$ be the row component $c_i$-adjacent to $r$. Without loss of generality, assume that $r$ and $s$ are arranged as follows:
\[ \begin{split}
r &= [a,b]_\Z \times \{k\}, \\
s &= [c,d]_\Z \times \{k+1\}.
\end{split}
\]
(Possibly $c$ or $d$ could be infinite.)

Let $P = [a,b]_\Z \cap [c,d]_\Z$, and let $Q = r - (P\times \{k\})$. This $Q$ has the structure of a path graph (perhaps of two components), and so the points of $Q$ can be sequentially deleted from $X$ by Lemma \ref{htp-equiv-remove-pt}, so $X$ is $2$-homotopy equivalent to $X-Q$.

If $i=2$, then the remaining points of $r$ (that is, the set $P\times \{k\}$) can all be removed by Lemma \ref{htp-equiv-remove-pt}, because we will have $N_{(X,c_2)}(y,k) \subseteq N_{(X,c_2)}(y,k+1)$ for each $y\in P$. Thus $X$ is $2$-homotopy equivalent to $X-r$, as desired.

If $i=1$, then the points of $P\times \{k\}$ can be removed by Lemma \ref{htp-equiv-remove-path}, where the paths $\gamma$ and $\delta$ in that lemma are parameterizations of $P\times \{k\}$ and $P\times \{k-1\}$.

Thus in either case, $(X,c_i)$ is $i$-homotopy equivalent to $(X-r,c_i)$, as desired.
\end{proof}


By inductively removing row components, we obtain:
\begin{lem}\label{tree-contractible}
Let $(X,c_i)$ be a digital image in $\Z^2$ such that the row graph $R(X,c_i)$ is a tree. Then $(X,c_i)$ is $i$-contractible.
\end{lem}
\begin{proof}
The proof is by induction on the number of row components of $(X,c_i)$. If there is only a single row component, then $X$ is isomorphic to some interval $[a,b]_\Z \times \{0\}$, which is $i$-contractible.

For the inductive step, assume the row component graph $R(X,c_i)$ is a tree of more than one vertex. Since it is a tree, there is some row component $r$ which is $c_i$-adjacent to only one other row component $s$. Then by Lemma \ref{htp-equiv-remove-row-component}, $(X,c_i)$ is $i$-homotopy equivalent to $(X-r,c_i)$, which is $i$-contractible by induction.
\end{proof}

\begin{thm}\label{noholecontractible}
Let $(X,c_i,c_j)$ be a digital picture. If $X$ is $c_i$-connected and has no hole, then $(X,c_i)$ is $i$-contractible.
\end{thm}
\begin{proof}
This follows immediately from Theorems \ref{rowgraphconnected} and \ref{rowgraphacyclic}, and Lemma \ref{tree-contractible}.
\end{proof}

\bigskip

\section{Row- and column-convex sets}\label{rcconvex}
\begin{definition}
For any $X\subset \Z^2$, the set $X$ is \emph{row-convex} when any horizontal row $X \cap (\Z\times \{k\})$ is connected. Similarly we say $X$ is \emph{column-convex} when any vertical column $X \cap (\{k\} \times \Z)$ is connected. When $X$ is both row-convex and column-convex, we call it \emph{row-column-convex}, or \emph{rc-convex}.

In some sources, e.g. \cite{bdnp01}, these concepts are called \emph{h-convex} and \emph{v-convex}. In the literature about polyominoes, e.g. \cite{dv84}, an rc-convex polyomino is simply called \emph{convex}. But the word ``convex'' has been used in the digital geometry and topology literature with different meanings (see \cite{kim82, boxe21} for two different uses of the word), so we use the more specific terms.
\end{definition}

\begin{thm}\label{rc-contractible}
Let $(X,c_i)$ be finite, $c_i$-connected, and row-convex. Then $(X,c_i)$ is $i$-contractible.
\end{thm}
\begin{proof}
Since $X$ is connected and row-convex, there is only one row component on each row, and thus the row component graph $R(X,c_i)$ is simply a path graph. Since this is a tree, $(X,c_i)$ is $i$-contractible by Lemma \ref{tree-contractible}.
\end{proof}

For a set $X\subset \Z^2$, let $H_r(X)$ be the \emph{row-convex hull of $X$}, defined as the intersection of all row-convex sets containing $X$. Equivalently, $H_r(X)$ is obtained from $X$ by ``filling in the gaps'' in any row, so that each row becomes connected. Clearly $X\subseteq H_r(X)$, and our next result will show that any digital picture $(X,c_i,c_j)$ with no hole is homotopy equivalent to $(H_r(X),c_i,c_j)$.

\begin{thm}\label{htp-equiv-row-hull}
Let $(X,c_i,c_j)$ be a digital picture such that $X$ is $c_i$-connected with no hole. Then $(X,c_i,c_j)$ is homotopy equivalent to $(H_r(X),c_i,c_j)$.
\end{thm}
\begin{proof}
By Lemma \ref{noholecontractible}, $(X,c_i)$ is $i$-contractible. Since $H_r(X)$ is row-convex, it is also $i$-contractible by Theorem \ref{rc-contractible}. Thus $(X,c_i)$ and $(H_r(X),c_i)$ are $i$-homotopy equivalent.

It remains to show that their complements are $j$-homotopy equivalent. Note that $\overline{H_r(X)}$ is the set obtained from $\bar X$ by removing all of the finite row components of $\bar X$. So it suffices to show that we can remove all finite row components of $\bar X$ by a $j$-homotopy equivalence.

We prove this by induction on the number of finite row components of $\bar X$. If there are no such finite components, then there is nothing to show.

For the inductive case, assume there is some finite row component. We will show that $(\bar X,c_j)$ is $j$-homotopy equivalent to $(\bar X-r,c_j)$ for some finite row component $r$ of $\bar X$.
Let $A$ be a maximal $c_j$-connected set of finite row components.
Because $X$ is connected, this set $A$ has no hole, and so Theorems \ref{rowgraphconnected} and \ref{rowgraphacyclic} show that the row component graph of $R(A,c_j)$ is a tree.

If $A$ includes more than one row component of $\bar X$, then because $R(A,c_j)$ is a tree we may choose a row component $r\subset \bar X$ which is $c_j$-adjacent to exactly one other row component of $\bar X$. Then by Lemma \ref{htp-equiv-remove-row-component}, $(\bar X,c_j)$ is $j$-homotopy equivalent to $(\bar X-r,c_j)$ as desired.

It remains to consider the case when $A$ consists only of a single row component which we call $q$. In this case we argue that $q$ is $c_j$-adjacent to exactly one other row-component of $\bar X$, and the same arguments above will complete the proof.

Since $X$ has no hole, this means $\bar X$ is $c_j$-connected, and thus $q \subset \bar X$ must be $c_j$-adjacent to some infinite row component of $\bar X$.

It remains to show that this $q$ cannot be $c_j$-adjacent to more than one infinite row component of $\bar X$. This can be verified by considering various cases that such a situation would induce, depending on the positioning of the two infinite row components. Assume that $q$ is $c_j$-adjacent to two distinct infinite row components $\bar X$ called $r_1,r_2$. If $r_1$ and $r_2$ are both ``above'' $q$, then we have a situation like in Figure \ref{r1r2qfig}. This situation is impossible because it will make $X$ disconnected, since the points of $X$ on the ends of $q$ cannot be joined by a path to the points of $X$ on the ends of $r_1$ and $r_2$. Similar arguments will show that any other configuration in which $q$ is $c_j$-adjacent to more than one infinite row component of $\bar X$ is impossible.
\begin{figure}
\[ \begin{tikzpicture}[scale=1]
\usetikzlibrary{arrows.meta}

\draw[{Bar}-{Bar}] (1.1,.5) -- (4.8,.5) node[pos=.5,fill=white] {$q$};
\draw[<-{Bar}] (-2,1.5) -- (1.8,1.5) node[pos=.5,fill=white] {$r_1$};
\draw[{Bar}->] (4.1,1.5) -- (8,1.5) node[pos=.5,fill=white] {$r_2$};
\foreach \x/\y in {-1/0,0/0,5/0,6/0,2/1,3/1} {
  \draw[fill=lightgray] (\x,\y) rectangle (\x+.9,\y+.9) {};
}
\end{tikzpicture}
\]
\caption{One case of an impossible configuration from the proof of Theorem \ref{htp-equiv-row-hull}. Shaded pixels are points of $X$.\label{r1r2qfig}}
\end{figure}
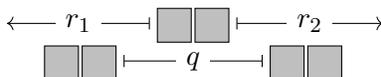
\end{proof}

All of the constructions presented above have focused on row components, row-convex sets, and the row-convex hull of sets. Repeating exactly the same arguments for columns instead of rows gives analogous results, in particular that any digital picture $(X,c_i,c_j)$ with no holes is homotopy equivalent to $(H_c(X),c_i,c_j)$, where $H_c(X)$ is the column-convex hull of $X$.

\begin{definition}
For a digital image $X\subset \Z^2$, we define $H_{rc}(X)$, the \emph{rc-convex hull} of $X$, as the intersection of all rc-convex sets containing $X$.
\end{definition}

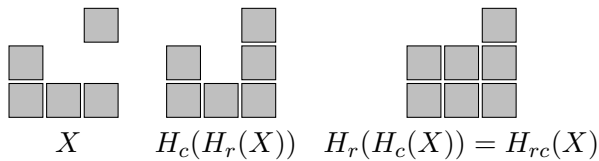
\begin{figure}
\[
\begin{tabular}{ccc}
\begin{tikzpicture}[scale=1]
\foreach \x/\y in {0/0,1/0,2/0,0/1,2/2} {
  \draw[fill=lightgray] (\x,\y) rectangle (\x+.9,\y+.9) {};
}
\end{tikzpicture} \qquad
&
\begin{tikzpicture}[scale=1]
\foreach \x/\y in {0/0,1/0,2/0,0/1,2/2,2/1} {
  \draw[fill=lightgray] (\x,\y) rectangle (\x+.9,\y+.9) {};
}
\end{tikzpicture} \qquad
&
\begin{tikzpicture}[scale=1]
\foreach \x/\y in {0/0,1/0,2/0,0/1,2/2,2/1,1/1} {
  \draw[fill=lightgray] (\x,\y) rectangle (\x+.9,\y+.9) {};
}
\end{tikzpicture}
\\
$X$ & $H_c(H_r(X))$ & $H_r(H_c(X)) = H_{rc}(X)$
\end{tabular}
\]
\caption{A digital image with $H_r(H_c(X)) \neq H_c(H_r(X))$.\label{rcrcfig}}
\end{figure}

The relationship of the rc-convex hull and the row- and column-convex hulls is not very straightforward. It is natural to expect that $H_{rc}(X) = H_r(H_c(X)) = H_c(H_r(X))$, but this is not generally true as shown in Figure \ref{rcrcfig}. When $X$ is $c_2$-connected, we are able to prove the following:
\begin{lem}\label{connectedrc}
If $X\subset \Z^2$ is $c_2$-connected, then $H_r(H_c(X)) = H_c(H_r(X))$.
\end{lem}
\begin{proof}
We will prove that $H_c(H_r(X)) \subseteq H_r(H_c(X))$. Symmetric arguments will show the other containment, which proves that these sets are equal. Take some $x\in H_c(H_r(X))$ and we will show $x\in H_r(H_c(X))$.

If $x\in H_c(X)$, then automatically $x\in H_c(H_r(X))$, so we may assume that $x\not\in H_c(X)$, so $x$ is an element of some ``row gap'' of $H_c(X)$. That is, we have some $a<b<c$ and $k$ with $x=(b,k)$ and $(a,k),(c,k) \in H_c(X)$.

If both of $(a,k)$ and $(c,k)$ are in $X$, then $x\in H_r(X) \subseteq H_r(H_c(X))$ as desired. So we assume without loss of generality that $(a,k) \not \in X$. In fact, for the same reason, we may assume $(d,k) \not \in X$ for all $d\le b$.

Note that $(a,k)\in H_c(X)$ but $(a,k)\not \in X$, and so $(a,k)$ is part of some ``column gap'' of $X$. That is, there are $j<k<l$ with $(a,j),(a,l)\in X\subseteq H_c(X)$.

\begin{figure}
\[
\begin{tikzpicture}[pixelmode]
\foreach \x/\y in {0/0,0/5} {
  \draw[fill=lightgray] (\x-.45,\y-.45) rectangle (\x+.45,\y+.45) {};
}
\node[below=.4cm] at (0,0) {$(a,j)\in X$};
\node[above=.4cm] at (0,5) {$(a,l)\in X$};

\foreach \x/\y in {5/3} {
  \draw[fill=white] (\x-.45,\y-.45) rectangle (\x+.45,\y+.45) {};
}
\node[right=.2cm] at (5,3) {$(c,k) \in H_c(X)$};

\foreach \x in {-1,...,2} {
 \node[comp-vertex] at (\x,3) {};
}
\node at (-2,3) {$\dots$};

\node[comp-vertex] at (2,3) {};
\node[above] at (2,3) {$(b,k)$};
\end{tikzpicture}
\]
\caption{Illustration for the proof of Lemma \ref{connectedrc}. All points marked with circles are not elements of $X$.\label{connectedrcfig}}
\end{figure}
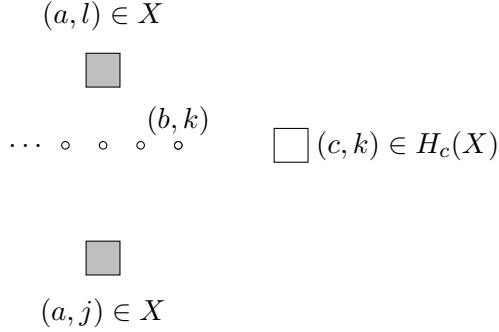

Summarizing above, we have $(a,j),(a,l) \in X$, with $j<l<k$ and $a<c$, but $(d,k)\not\in X$ for all $d\le b < c$. The situation resembles Figure \ref{connectedrcfig}. Our assumption that $X$ is $c_2$-connected means that there is a $c_2$-path in $X$ from $(a,l)$ to $(a,j)$, and this path must necessarily pass through some points $(b,m)$ and $(b,n)$ with $n<l<m$. Thus we have $x=(b,k) \in H_c(X) \subseteq H_r(H_c(X))$ as desired.
\end{proof}

We still would like to understand the relationship between $H_{rc}(X)$ and the sets $H_r(H_c(X))$ and $H_c(H_r(X))$. Given a set $X\subset \Z^2$, we may apply the row- and column-convex hulls $H_r$ and $H_c$ iteratively to obtain nested sets:
\begin{equation}\label{rccontainments}
X \subseteq H_r(X) \subseteq H_c(H_r(X)) \subseteq H_r(H_c(H_r(X))) \subseteq \dots
\end{equation}
When $X$ is finite, all of these sets above will be contained in the ``bounding box'' of $X$, the smallest rectangle containing $X$. Since these sets are all contained in the finite bounding box, the containments of \eqref{rccontainments} must eventually stabilize. We let $H(X)$ denote the set to which the containments of \eqref{rccontainments} stabilizes.

\begin{figure}
\tikzset{every picture/.append style={scale=.4}}
\begin{align*}
\vcbox{
\begin{tikzpicture}[pixelmode]
\foreach \x/\y in {0/0,0/1,2/3,4/0,5/5,6/2,7/7,9/4,11/6} {
  \draw[fill=lightgray] (\x,\y) rectangle (\x+.9,\y+.9) {};
}
\end{tikzpicture}}
&\overset{H_r}\longrightarrow
\vcbox{
\begin{tikzpicture}[pixelmode]
\foreach \x/\y in {0/0,1/0,2/0,3/0,4/0,0/1,6/2,2/3,9/4,5/5,11/6,7/7} {
  \draw[fill=lightgray] (\x,\y) rectangle (\x+.9,\y+.9) {};
}
\end{tikzpicture}}
\overset{H_c}\longrightarrow
\vcbox{
\begin{tikzpicture}[pixelmode]
\foreach \x/\y in {0/0,0/1,1/0,2/0,2/1,2/2,2/3,3/0,4/0,5/5,6/2,7/7,9/4,11/6} {
  \draw[fill=lightgray] (\x,\y) rectangle (\x+.9,\y+.9) {};
}
\end{tikzpicture}}
\overset{H_r}\longrightarrow
\vcbox{
\begin{tikzpicture}[pixelmode]
\foreach \x/\y in {0/0,1/0,2/0,3/0,4/0,0/1,1/1,2/1,2/2,3/2,4/2,5/2,6/2,2/3,9/4,5/5,11/6,7/7} {
  \draw[fill=lightgray] (\x,\y) rectangle (\x+.9,\y+.9) {};
}
\end{tikzpicture}}
\\
&\overset{H_c}\longrightarrow
\vcbox{
\begin{tikzpicture}[pixelmode]
\foreach \x/\y in {0/0,0/1,1/0,1/1,2/0,2/1,2/2,2/3,3/0,3/1,3/2,4/0,4/1,4/2,5/2,5/3,5/4,5/5,6/2,7/7,9/4,11/6} {
  \draw[fill=lightgray] (\x,\y) rectangle (\x+.9,\y+.9) {};
}
\end{tikzpicture}}
\overset{H_r}\longrightarrow
\vcbox{
\begin{tikzpicture}[pixelmode]
\foreach \x/\y in {0/0,1/0,2/0,3/0,4/0,0/1,1/1,2/1,3/1,4/1,2/2,3/2,4/2,5/2,6/2,2/3,3/3,4/3,5/3,5/4,6/4,7/4,8/4,9/4,5/5,11/6,7/7} {
  \draw[fill=lightgray] (\x,\y) rectangle (\x+.9,\y+.9) {};
}
\end{tikzpicture}
}
\overset{H_c}\longrightarrow
\vcbox{
\begin{tikzpicture}[pixelmode]
\foreach \x/\y in {0/0,0/1,1/0,1/1,2/0,2/1,2/2,2/3,3/0,3/1,3/2,3/3,4/0,4/1,4/2,4/3,5/2,5/3,5/4,5/5,6/2,6/3,6/4,7/4,7/5,7/6,7/7,8/4,9/4,11/6} {
  \draw[fill=lightgray] (\x,\y) rectangle (\x+.9,\y+.9) {};
}
\end{tikzpicture}
}
\\
&\overset{H_r}\longrightarrow
\vcbox{
\begin{tikzpicture}[pixelmode]
\foreach \x/\y in {0/0,1/0,2/0,3/0,4/0,0/1,1/1,2/1,3/1,4/1,2/2,3/2,4/2,5/2,6/2,2/3,3/3,4/3,5/3,6/3,5/4,6/4,7/4,8/4,9/4,5/5,6/5,7/5,7/6,8/6,9/6,10/6,11/6,7/7} {
  \draw[fill=lightgray] (\x,\y) rectangle (\x+.9,\y+.9) {};
}
\end{tikzpicture}
}
\overset{H_c}\longrightarrow
\vcbox{
\begin{tikzpicture}[pixelmode]
\foreach \x/\y in {0/0,0/1,1/0,1/1,2/0,2/1,2/2,2/3,3/0,3/1,3/2,3/3,4/0,4/1,4/2,4/3,5/2,5/3,5/4,5/5,6/2,6/3,6/4,6/5,7/4,7/5,7/6,7/7,8/4,8/5,8/6,9/4,9/5,9/6,10/6,11/6} {
  \draw[fill=lightgray] (\x,\y) rectangle (\x+.9,\y+.9) {};
}
\end{tikzpicture}
}
\end{align*}
\caption{The set $X$ shown above requires many iterated applications of $H_r$ and $H_c$ before stabilizing to the set $H(X)$. Extending the structure of $X$ allows for arbitrarily long iterations of $H_r$ and $H_c$.\label{longrcfig}}
\end{figure}
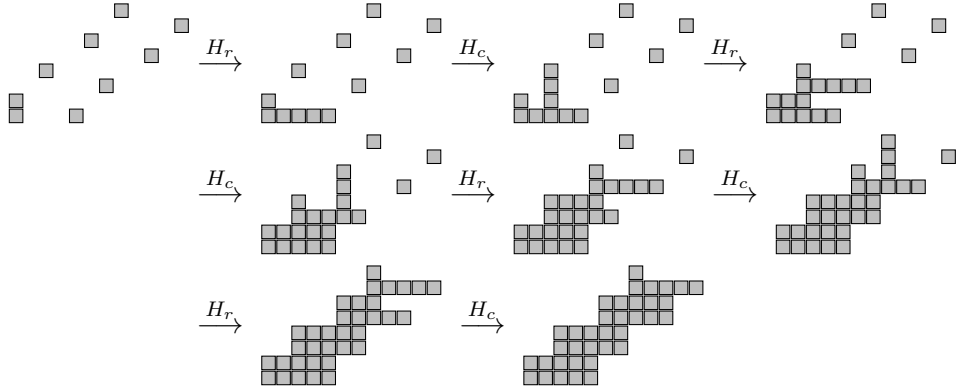

As shown in Figure \ref{longrcfig}, this sequence of containments can be arbitrarily long before it stabilizes. Due to Lemma \ref{connectedrc} this kind of example cannot happen when $X$ is connected.

Now we show that the set $H(X)$ in fact equals the rc-convex hull $H_{rc}(X)$.
\begin{thm}\label{H=Hrc}
Let $X\subset \Z^2$ be a finite set. Then $H(X) = H_{rc}(X)$.
\end{thm}
\begin{proof}
Since $H(X)$ is stable under application of $H_r$ and $H_c$, this means that $H(X)$ is rc-convex and thus we have $H_{rc}(X) \subseteq H(X)$. So we must show $H(X) \subseteq H_{rc}(X)$. We demonstrate this by showing that all of the sets in the nested sequence \eqref{rccontainments} are contained in $H_{rc}(X)$.

Let $S_k$ be the $k$th set in the nested sequence \eqref{rccontainments}, so that $S_0 = X$, $S_1 = H_r(X)$, etc. We will show that $S_k\subseteq H_{rc}(X)$ by induction on $k$. The case $k=0$ is clear because $X\subseteq H_{rc}(X)$.

For the inductive case, assume without loss of generality that $k$ is odd, so that $S_k = H_r(S_{k-1})$. Similar arguments will hold when $k$ is even and $S_k = H_c(S_{k-1})$. By induction we may assume that $S_{k-1}\subseteq H_{rc}(X)$. To show $S_k \subseteq H_{rc}(X)$, take $x\in S_k$ and we will show $x\in H_{rc}(X)$.

If $x\in S_{k-1}$ then we have $x\in H_{rc}(X)$, so $x$ is an element of every rc-convex set containing $X$. So we assume that $x\not \in S_{k-1}$. This means that $x$ is an element of some ``row-gap'' of $S_{k-1}$. That is, there are integers $a<b<c$ and $l$ with $x=(b,l)$ and $(a,l),(c,l) \in S_{k-1} \subset H_{rc}(X)$. Since $H_{rc}(X)$ is row-convex, this means that $x \in H_{rc}(X)$ as desired.
\end{proof}

When $X\subset \Z^2$ is $c_2$-connected, we obtain the expected result:

\begin{thm}
If $X\subset \Z^2$ is $c_2$-connected, then $H_{rc}(X) = H_r(H_c(X)) = H_c(H_r(X))$.
\end{thm}
\begin{proof}
Lemma \ref{connectedrc} has already shown that $H_r(H_c(X)) = H_c(H_r(X))$. To finish the proof, note that we will have:
\[ H_r(H_c(H_r(X))) = H_r(H_r(H_c(X))) = H_r(H_c(X)) = H_c(H_r(X)) \]
because $H_r$ is idempotent. Thus $H_c(H_r(X))$ is the stable set $H(X)$ in the nested containments \eqref{rccontainments}. Thus by Theorem \ref{H=Hrc} we have $H_c(H_r(X)) = H_{rc}(X)$ as desired.
\end{proof}

The results above will combine to give:
\begin{thm}\label{htp-equiv-Hrc}
Any $c_i$-connected digital picture $(X,c_i,c_j)$ with no hole is homotopy equivalent to $(H_{rc}(X),c_i,c_j)$.
\end{thm}
\begin{proof}
By Theorem \ref{htp-equiv-row-hull} and the corresponding statement for column-convex hulls, taking the row- or column-convex hull of any digital picture does not change the homotopy type. Thus all sets in the sequence \eqref{rccontainments} are homotopy equivalent to each other, and thus $(X,c_i,c_j)$ is homotopy equivalent to $(H(X),c_i,c_j)$, which by Theorem \ref{H=Hrc} equals $(H_{rc}(X),c_i,c_j)$.
\end{proof}

\bigskip

\section{The homotopy type of a digital picture}\label{homotopytypeofdigitalpicture}
The paper \cite{hmps15} discussed ``reducible'' digital images. A finite digital image $(X,\kappa)$ is $i$-reducible when it is $i$-homotopy equivalent to a digital image of fewer points. The following result appears in \cite{hmps15}:
\begin{lem}[\cite{hmps15}, Lemma 2.8]\label{reduciblenonsurj}
A finite digital image $(X,\kappa)$ is $i$-reducible if and only if $\id_X$ is $i$-homotopic in one step to a non-surjection.
\end{lem}
(The paper \cite{hmps15} deals exclusively with $1$-homotopy, but the proof given for the result above holds as-written for both $i\in \{1,2\}$.)

\begin{thm}
Let $X$ be a finite digital image, and $i\in \{1,2\}$. Then $X$ is $i$-reducible if and only if it is $i$-homotopy equivalent to a proper subset of itself.
\end{thm}
\begin{proof}
Certainly if $X$ is $i$-homotopy equivalent to a proper subset of itself, then it is $i$-homotopy equivalent to a set of fewer points and so is $i$-reducible.

For the other implication, assume that $X$ is $i$-reducible. Then by Lemma \ref{reduciblenonsurj} there is a nonsurjection $f:X\to X$ with $f\simeq_i \id_X$. Then by Lemma \ref{htp-equiv-image}, $X$ is $i$-homotopy equivalent to $f(X)$ which is a proper subset of $X$.
\end{proof}

For $2$-homotopy we can be more specific:
\begin{thm}\label{2reducible1point}
Let $X$ be a finite digital image. If $X$ is $2$-reducible, then there is some $p\in X$ such that  $X - \{p\}$ is $2$-homotopy equivalent to $X$.
\end{thm}
\begin{proof}
Assume that $X$ is 2-reducible. Then by Lemma \ref{reduciblenonsurj}, there is a nonsurjection $f:X\to X$ with $\id_X \simeq_2 f$. By the Spider Move Lemma we can take this to be a homotopy through spider moves. Let $g$ be the first non-identity stage of the homotopy, and since $g$ is a spider move we will have $g(X) = X-\{p\}$ for some $p$. Then by Lemma \ref{htp-equiv-image}, $X-\{p\}$ is 2-homotopy equivalent to $X$.
\end{proof}

We do not know if the same statement is true for 1-reducibility:
\begin{quest}
Let $X$ be a $1$-reducible digital image. Then is there some $p\in X$ such that $X$ is $1$-homotopy equivalent to $X-\{p\}$?
\end{quest}

The obvious interesting example when considering the question above is the square $X=[0,1]_\Z\times [0,1]_\Z$ with $c_1$-adjacency, which is 1-contractible. It can be checked that there is no selfmap $f:X\to X$ such that $f(X)$ has cardinality 3. In this case $X$ is indeed homotopy equivalent to $X-\{(1,1)\}$ (because both are 1-contractible), but not in a single step.

The discussion in \cite{hmps15} focused on the possibility of reducing the cardinality of a digital image by homotopy equivalence. \emph{Increasing} the cardinality was never discussed, because it is always possible to add points to a digital image (by introducing tree-like spurs, for example) without changing its homotopy type. The same is not true, however, for digital pictures.

\begin{thm}
A given digital picture $(X,c_i,c_j)$ is homotopy equivalent to only finitely many other digital pictures. Thus, among all digital pictures homotopy equivalent to $(X,c_i,c_j)$, there is one (possibly nonunique) of minimal cardinality, and one (possibly not unique) of maximal cardinality.
\end{thm}
\begin{proof}
Recall that if $(X,c_i,c_j)$ is homotopy equivalent to $(Y,c_i,c_j)$, then we must have $O(X,c_i,c_j) = O(Y,c_i,c_j)$. But this invariant represents the minimal length of a loop that surrounds $X$, and there are only finitely many digital images which can be surrounded of a loop of this length. Thus there are only finitely many such digital pictures $(Y,c_i,c_j)$.
\end{proof}

The smallest digital picture in the homotopy class is typically larger than the smallest digital image in the homotopy class. For example the interval $(I_n,c_i)$ is $2$-contractible and so the smallest digital image in its $2$-homotopy class is a single point. But $(I_n,c_i,c_j)$ is not homotopy equivalent to any proper subset, since any proper subset is either disconnected, or has lesser outer perimeter.

We would also like to understand the largest digital picture in the homotopy class. Theorem \ref{htp-equiv-Hrc} shows that a digital picture without holes can always be enlarged to its rc-convex hull. But in some cases the image can be enlarged even further:
\begin{exa}
The digital picture $(I_{3,4},c_1,c_2)$ equals its own rc-convex hull, but it is homotopy equivalent to the larger digital picture $(X,c_1,c_2)$ shown in Figure \ref{enlarge-c1-fig}. The 1-homotopy equivalence of $I_{3,4}$ and $X$ can be achieved by repeated application of Lemmas \ref{htp-equiv-remove-pt} and \ref{htp-equiv-remove-path}. The 2-homotopy equivalence of $\bar I_{3,4}$ and $\bar X$ can be achieved by repeated application of Lemma \ref{htp-equiv-remove-pt}.

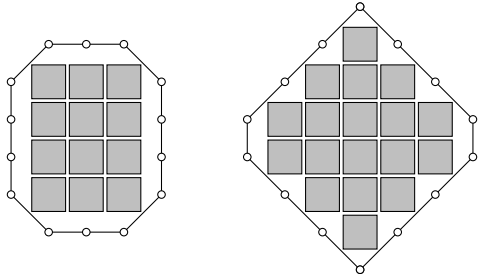
\begin{figure}
\[
\vcbox{
\begin{tikzpicture}[pixelmode]
\foreach \x/\y in {0/0,0/1,0/2,0/3,1/0,1/1,1/2,1/3,2/0,2/1,2/2,2/3} {
  \draw[fill=lightgray] (\x,\y) rectangle (\x+.9,\y+.9) {};
}
\begin{scope}[shift={(.45,.45)}]
\draw[edge] (-1, 3)--(0, 4)--(1, 4)--(2, 4)--(3, 3)--(3, 2)--(3, 1)--(3, 0)--(2, -1)--(1, -1)--(0, -1)--(-1, 0)--(-1, 1)--(-1, 2)--(-1, 3);
\foreach \x/\y in {-1/0,-1/1,-1/2,-1/3,3/0,3/1,3/2,3/3,0/-1,0/4,1/-1,1/4,2/-1,2/4} {
  \node[comp-vertex] at (\x,\y) {};
}
\end{scope}
\end{tikzpicture}
}
\qquad
\vcbox{
\begin{tikzpicture}[pixelmode]
\foreach \x/\y in {0/0,0/1,0/2,0/3,1/0,1/1,1/2,1/3,2/0,2/1,2/2,2/3,1/-1,1/4,-1/1,-1/2,3/1,3/2} {
  \draw[fill=lightgray] (\x,\y) rectangle (\x+.9,\y+.9) {};
}
\begin{scope}[shift={(.45,.45)}]
\draw[edge] (0, 4)--(1, 5)--(2, 4)--(3, 3)--(4, 2)--(4, 1)--(3, 0)--(2, -1)--(1, -2)--(0, -1)--(-1, 0)--(-2, 1)--(-2, 2)--(-1, 3)--(0, 4);
\foreach \x/\y in {-2/2,-1/3,0/4,1/5,1/5,2/4,3/3,4/2,-2/1,-1/0,0/-1,1/-2,1/-2,2/-1,3/0,4/1} {
  \node[comp-vertex] at (\x,\y) {};
}
\end{scope}
\end{tikzpicture}
}
\]
\caption{$(I_{3,4},c_1,c_2)$ at left is rc-convex, but homotopy equivalent to the larger digital picture $(X,c_1,c_2)$ at right. Each time we show a loop in the complement realizing the outer perimeter.\label{enlarge-c1-fig}}
\end{figure}
\end{exa}

In the example above the diagonal adjacencies in the complement allowed $I_{3,4}$ to be enlarged without increasing the outer perimeter. We do not believe this to be possible in a digital picture of the form $(X,c_2,c_1)$, so we make the following conjecture:
\begin{conj}
If $(X,c_2,c_1)$ is a digital picture with no hole, then the rc-convex hull $(H_{rc}(X),c_2,c_1)$ is a set of largest cardinality in the homotopy type of $(X,c_2,c_1)$.
\end{conj}

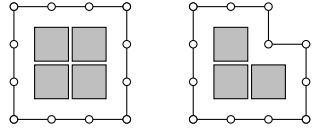
\begin{figure}
\[ 
\begin{tikzpicture}[pixelmode]
\foreach \x/\y in {0/0,0/1,1/0,1/1} {
  \draw[fill=lightgray] (\x,\y) rectangle (\x+.9,\y+.9) {};
}
\begin{scope}[shift={(.45,.45)}]
\draw[edge] (-1, -1)--(0, -1)--(1, -1)--(2, -1)--(2, 0)--(2, 1)--(2, 2)--(1, 2)--(0, 2)--(-1, 2)--(-1, 1)--(-1, 0)--(-1, -1);
\foreach \x/\y in {-1/-1,-1/2,0/-1,0/2,1/-1,1/2,2/-1,2/2,-1/-1,-1/0,-1/1,-1/2,2/-1,2/0,2/1,2/2} {
  \node[comp-vertex] at (\x,\y) {};
}
\end{scope}
\end{tikzpicture}
\qquad
\begin{tikzpicture}[pixelmode]
\foreach \x/\y in {0/0,0/1,1/0} {
  \draw[fill=lightgray] (\x,\y) rectangle (\x+.9,\y+.9) {};
}
\begin{scope}[shift={(.45,.45)}]
\draw[edge] (-1, -1)--(0, -1)--(1, -1)--(2, -1)--(2, 0)--(2, 1)--(1,1)--(1, 2)--(0, 2)--(-1, 2)--(-1, 1)--(-1, 0)--(-1, -1);
\foreach \x/\y in {-1/-1,-1/2,0/-1,0/2,1/-1,1/2,2/-1,1/1,-1/-1,-1/0,-1/1,-1/2,2/-1,2/0,2/1} {
  \node[comp-vertex] at (\x,\y) {};
}
\end{scope}
\end{tikzpicture}
\]
\caption{These digital pictures have the same outer perimeter, but do not appear to be homotopy equivalent.\label{I22fig}}
\end{figure}

We remark that the outer perimeter is not the only obstacle to enlarging a digital picture by a homotopy equivalence. Consider for example the digital pictures $(I_{2,2},c_2,c_1)$ and $(X,c_2,c_1)$ where $X = I_{2,2} - \{(2,2)\}$, shown in Figure \ref{I22fig}. Here we have $O(I_{2,2},c_2,c_1) = O(X,c_2,c_1) = 12$, and $(X,c_2)$ and $(I_{2,2},c_2)$ are $2$-homotopy equivalent. But $(\overline{I_{2,2}},c_1)$ and $(\bar X,c_1)$ do not seem to be 1-homotopy equivalent.

Thus a natural question arises:
\begin{quest}
Given a digital picture $(X,c_i,c_j)$, find some smallest and largest digital pictures $(Y,c_i,c_j)$ which are homotopy equivalent to $(X,c_i,c_j)$.
\end{quest}

\begin{figure}
\captionsetup[subfigure]{labelformat=empty}
{\tikzset{every picture/.append style={scale=.6}}
\begin{subfigure}[t]{.3\textwidth}
\begin{tikzpicture}
\foreach \x/\y in {0/0,7/0,8/0,9/0,0/1,1/1,2/1,5/1,7/1,9/1,2/2,3/2,4/2,5/2,6/2,7/2,8/2,9/2,5/3,8/3,8/4} {
  \draw[fill=lightgray] (\x,\y) rectangle (\x+.9,\y+.9) {};
}
\begin{scope}[shift={(.45,.45)}]
\draw[edge] (5, 4)--(6, 4)--(7, 4)--(8, 5)--(9, 4)--(10, 3)--(10, 2)--(10, 1)--(10, 0)--(9, -1)--(8, -1)--(7, -1)--(6, -1)--(5, -1)--(4, -1)--(3, -1)--(2, -1)--(1, -1)--(0, -1)--(-1, 0)--(-1, 1)--(0, 2)--(1, 2)--(2, 3)--(3, 3)--(4, 3)--(5, 4);
\foreach \x/\y in {-1/0,-1/1,0/2,1/2,2/3,3/3,4/3,5/4,6/4,7/4,8/5,9/4,10/3,10/2,10/1,10/0,0/-1,1/-1,2/-1,3/-1,4/-1,5/-1,6/-1,7/-1,8/-1,9/-1} {
  \node[comp-vertex] at (\x,\y) {};
}
\end{scope}
\path (0,-5) -- (5,5);
\end{tikzpicture}
\subcaption{$(X,c_1,c_2)$}
\end{subfigure}
\hfill
\begin{subfigure}[t]{.3\textwidth}
\begin{tikzpicture}
\foreach \x/\y in {0/0,7/0,8/0,9/0,0/1,1/1,2/1,7/1,9/1,2/2,3/2,4/2,5/2,6/2,7/2,8/2,9/2,8/3,8/4} {
  \draw[fill=lightgray] (\x,\y) rectangle (\x+.9,\y+.9) {};
}
\begin{scope}[shift={(.45,.45)}]
\draw[edge] (5, 3)--(6, 3)--(7, 4)--(8, 5)--(9, 4)--(10, 3)--(10, 2)--(10, 1)--(10, 0)--(9, -1)--(8, -1)--(7, -1)--(6, -1)--(5, -1)--(4, -1)--(3, -1)--(2, -1)--(1, -1)--(0, -1)--(-1, 0)--(-1, 1)--(0, 2)--(1, 2)--(2, 3)--(3, 3)--(4, 3)--(5, 3);
\foreach \x/\y in {-1/0,-1/1,0/2,1/2,2/3,3/3,4/3,5/3,6/3,7/4,8/5,9/4,10/3,10/2,10/1,10/0,0/-1,1/-1,2/-1,3/-1,4/-1,5/-1,6/-1,7/-1,8/-1,9/-1} {
  \node[comp-vertex] at (\x,\y) {};
}
\end{scope}
\path (0,-5) -- (5,5);
\end{tikzpicture}
\subcaption{Minimal reduction of $(X,c_1,c_2)$}
\end{subfigure}
\hfill
\begin{subfigure}[t]{.3\textwidth}
\begin{tikzpicture}[pixelmode]
\foreach \x/\y in {0/0,1/0,2/0,3/0,4/0,5/0,6/0,7/0,8/0,9/0,1/-1,2/-1,3/-1,4/-1,5/-1,6/-1,7/-1,8/-1,2/-2,3/-2,4/-2,5/-2,6/-2,7/-2,3/-3,4/-3,5/-3,6/-3,4/-4,5/-4,0/1,1/1,2/1,3/1,4/1,5/1,6/1,7/1,9/1,10/1,1/2,2/2,3/2,4/2,5/2,6/2,7/2,8/2,9/2,10/2,2/3,3/3,4/3,5/3,6/3,7/3,8/3,9/3,3/4,4/4,5/4,6/4,7/4,8/4,4/5,5/5,6/5,7/5,5/6,6/6} {
  \draw[fill=lightgray] (\x,\y) rectangle (\x+.9,\y+.9) {};
}
\begin{scope}[shift={(.45,.45)}]
\draw[edge] (5, -5)--(6, -4)--(7, -3)--(8, -2)--(9, -1)--(10, 0)--(11, 1)--(11, 2)--(10, 3)--(9, 4)--(8, 5)--(7, 6)--(6, 7)--(5, 7)--(4, 6)--(3, 5)--(2, 4)--(1, 3)--(0, 2)--(-1, 1)--(-1, 0)--(0, -1)--(1, -2)--(2, -3)--(3, -4)--(4, -5)--(5, -5);
\foreach \x/\y in {-1/0,0/-1,1/-2,2/-3,3/-4,4/-5,5/-5,6/-4,7/-3,8/-2,9/-1,10/0,11/1,11/2,10/3,9/4,8/5,7/6,6/7,-1/1,0/2,1/3,2/4,3/5,4/6,5/7} {
  \node[comp-vertex] at (\x,\y) {};
}
\end{scope}
\end{tikzpicture}
\subcaption{Maximal expansion of $(X,c_1,c_2)$}
\end{subfigure}
\begin{subfigure}[b]{.3\textwidth}
\begin{tikzpicture}
\foreach \x/\y in {0/0,7/0,8/0,9/0,0/1,1/1,2/1,5/1,7/1,9/1,2/2,3/2,4/2,5/2,6/2,7/2,8/2,9/2,5/3,8/3,8/4} {
  \draw[fill=lightgray] (\x,\y) rectangle (\x+.9,\y+.9) {};
}
\begin{scope}[shift={(.45,.45)}]
\draw[edge] (2, -1)--(1, -1)--(0, -1)--(-1, -1)--(-1, 0)--(-1, 1)--(-1, 2)--(0, 2)--(1, 2)--(1, 3)--(2, 3)--(3, 3)--(4, 3)--(4, 4)--(5, 4)--(6, 4)--(7, 4)--(7, 5)--(8, 5)--(9, 5)--(9, 4)--(9, 3)--(10, 3)--(10, 2)--(10, 1)--(10, 0)--(10, -1)--(9, -1)--(8, -1)--(7, -1)--(6, -1)--(5, -1)--(4, -1)--(3, -1)--(2, -1);
\foreach \x/\y in {-1/0,-1/1,-1/2,0/2,1/2,1/3,2/3,3/3,4/3,4/4,5/4,6/4,7/4,7/5,8/5,9/5,9/4,9/3,10/3,10/2,10/1,10/0,10/-1,-1/-1,0/-1,1/-1,2/-1,3/-1,4/-1,5/-1,6/-1,7/-1,8/-1,9/-1} {
  \node[comp-vertex] at (\x,\y) {};
}
\end{scope}
\end{tikzpicture}
\subcaption{$(X,c_2,c_1)$}
\end{subfigure}
\hfill
\begin{subfigure}[b]{.3\textwidth}
\begin{tikzpicture}
\foreach \x/\y in {0/0,7/0,8/0,9/0,0/1,1/1,2/1,7/1,9/1,2/2,3/2,4/2,5/2,6/2,7/2,8/2,9/2,5/3,8/3,8/4} {
  \draw[fill=lightgray] (\x,\y) rectangle (\x+.9,\y+.9) {};
}
\begin{scope}[shift={(.45,.45)}]
\draw[edge] (2, -1)--(1, -1)--(0, -1)--(-1, -1)--(-1, 0)--(-1, 1)--(-1, 2)--(0, 2)--(1, 2)--(1, 3)--(2, 3)--(3, 3)--(4, 3)--(4, 4)--(5, 4)--(6, 4)--(7, 4)--(7, 5)--(8, 5)--(9, 5)--(9, 4)--(9, 3)--(10, 3)--(10, 2)--(10, 1)--(10, 0)--(10, -1)--(9, -1)--(8, -1)--(7, -1)--(6, -1)--(5, -1)--(4, -1)--(3, -1)--(2, -1);
\foreach \x/\y in {-1/0,-1/1,-1/2,0/2,1/2,1/3,2/3,3/3,4/3,4/4,5/4,6/4,7/4,7/5,8/5,9/5,9/4,9/3,10/3,10/2,10/1,10/0,10/-1,-1/-1,0/-1,1/-1,2/-1,3/-1,4/-1,5/-1,6/-1,7/-1,8/-1,9/-1} {
  \node[comp-vertex] at (\x,\y) {};
}
\end{scope}
\end{tikzpicture}
\subcaption{Minimal reduction of $(X,c_2,c_1)$}
\end{subfigure}
\hfill
\begin{subfigure}[b]{.3\textwidth}
\begin{tikzpicture}[pixelmode]
\foreach \x/\y in {0/0,1/0,2/0,3/0,4/0,5/0,6/0,7/0,8/0,9/0,0/1,1/1,2/1,3/1,4/1,5/1,6/1,7/1,9/1,2/2,3/2,4/2,5/2,6/2,7/2,8/2,9/2,5/3,6/3,7/3,8/3,8/4} {
  \draw[fill=lightgray] (\x,\y) rectangle (\x+.9,\y+.9) {};
}
\begin{scope}[shift={(.45,.45)}]
\draw[edge] (2, -1)--(1, -1)--(0, -1)--(-1, -1)--(-1, 0)--(-1, 1)--(-1, 2)--(0, 2)--(1, 2)--(1, 3)--(2, 3)--(3, 3)--(4, 3)--(4, 4)--(5, 4)--(6, 4)--(7, 4)--(7, 5)--(8, 5)--(9, 5)--(9, 4)--(9, 3)--(10, 3)--(10, 2)--(10, 1)--(10, 0)--(10, -1)--(9, -1)--(8, -1)--(7, -1)--(6, -1)--(5, -1)--(4, -1)--(3, -1)--(2, -1);
\foreach \x/\y in {-1/0,-1/1,-1/2,0/2,1/2,1/3,2/3,3/3,4/3,4/4,5/4,6/4,7/4,7/5,8/5,9/5,9/4,9/3,10/3,10/2,10/1,10/0,10/-1,-1/-1,0/-1,1/-1,2/-1,3/-1,4/-1,5/-1,6/-1,7/-1,8/-1,9/-1} {
  \node[comp-vertex] at (\x,\y) {};
}
\end{scope}
\end{tikzpicture}
\subcaption{Maximal expansion of $(X,c_2,c_1)$}
\end{subfigure}
}
\caption{A set $X\subset \Z^2$ and the apparent smallest and largest digital pictures in its homotopy type. Each time we show a loop in $\bar X$ realizing the outer perimeter.\label{biggestsmallestfig}}
\end{figure}
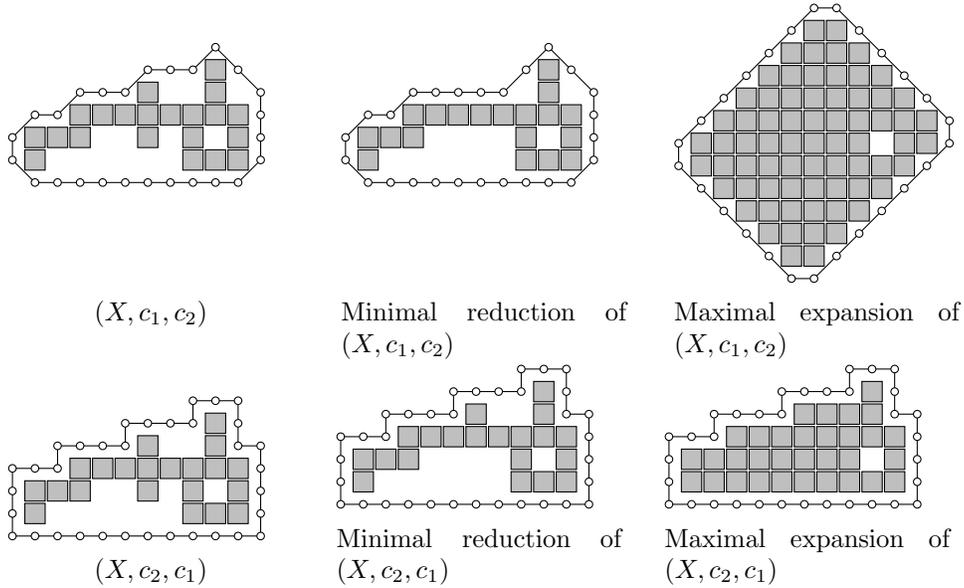

Figure \ref{biggestsmallestfig} shows an attempt to do this for the digital picture from Figure \ref{rowgraphfig}, though the results are only conjectural. We have not proven that they are in fact the smallest and largest digital pictures in the homotopy class: we present the figure only to give an impression of the problem.

\bigskip

\bibliography{digital}
\bibliographystyle{plain}

\end{document}